\documentclass[12pt,a4paper]{article}

\usepackage{amsfonts,amssymb,wasysym,wrapfig}
\usepackage{graphicx}

\newcommand{\Co}{\mathbb{C}}

\newcommand{\R}{\mathbb{R}}

\newtheorem{theorem}{Theorem}[section]

\newtheorem{lemma}[theorem]{Lemma}
\newtheorem{proposition}[theorem]{Proposition}

\newtheorem{remark}[theorem]{Remark}

\newtheorem{conjecture}[theorem]{Conjecture}

\title{How to Fold a Manifold}

\author{J. Scott Carter\\ University of South Alabama \and Seiichi Kamada \\ Hiroshima University}

\begin{document}
\maketitle
 
 \section{Introduction}
The English word ``manifold''  evokes images of an object that is layered or folded. Of course, the mathematical definition is quite different where the emphasis is upon ``many" such as the multiplicity of descriptions afforded by coordinate charts or indeed the need to define many coordinate systems --- one for each point in the space. Nevertheless, there are structures in which the colloquial meaning and the mathematical meaning coincide. Specifically, in the case of a covering space or a branched cover, we imagine that the manifold is stacked in layers above the base space. According to a classical theorem of Alexander \cite{Al1920}, every compact connected orientable manifold can be expressed as an irregular simple branched covering of the sphere with branch loci consisting of a codimension $2$ subcomplex. In dimensions two through four, we assume that the branch locus is  a finite collection of points, a knot or link, or a knotted or linked surface \cite{Hi1976,Mo1976,IP2002}. 
In this paper, we demonstrate methods of folding these coverings. We prove the following results.

\begin{theorem}\label{twofold} Let $k=-1,0,1,$ or $2$. Let $f: M^{k+2} \rightarrow S^{k+2}$ be a $2$-fold branched cover of $S^{k+2}$ branched along a closed oriented (possibly disconnected) smooth sub-manifold $L^k \hookrightarrow S^{k+2}$. Then there is an embedding $\tilde{f}: M^{k+2}\hookrightarrow S^{k+2} \times D^2$ such that $p \circ \tilde{f} = f$, where $p :  S^{k+2} \times D^2  \rightarrow S^{k+2}$ is the projection onto the first factor. 
 \end{theorem}
 
 Such an embedding, $\tilde{f}$, will be called an {\it folded embedding} or an {\it embedded folding}. We will say that the cover is {\it folded}. In case of $n$-fold irregular simple branched coverings $f: M^{k+2} \rightarrow S^{k+2}$ , there are topological obstructions to achieving an embedding.   Still in many cases, we can construct an immersion $\tilde{f}: M^{k+2} \looparrowright S^{k+2}\times D^2$ such that the composition $p\circ \tilde{f}=f$. In this case, we call $\tilde{f}$ an {\it immersed folding} or a {\it folded immersion}.

\begin{theorem}
\label{mainimmersed} Let $k=-1,0,$ or $1.$ Let $f: M^{k+2} \rightarrow S^{k+2}$ be a $3$-fold simple branched cover of $S^{k+2}$  branched along a closed  oriented (possibly disconnected) smooth sub-manifold $L^k \hookrightarrow S^{k+2}$. Then there is an immersion $\tilde{f}: M^{k+2} \looparrowright S^{k+2} \times D^2$ such that the restriction of the projection onto the first factor is the covering map $f$. Thus the lift $\tilde{f}$ is an immersed folding.
\end{theorem}

Often our techniques work for branched coverings of degree greater than $3$ and in case $k=2$. In this paper, we concentrate on the lower degree branched covers. We postpone the proof of Theorem~\ref{mainimmersed} in case $k=2$ for a subsequent paper.  We also give examples of foldings of $3$ and $4$ dimensional spheres that are constructed from our techniques. The advantages to our constructions are that the embeddings and immersions are given by explicit descriptions for which invariants such as fundamental group or Fox colorings can easily be computed. 

To illustrate the problems associated with constructing embedded foldings, consider 
the knot $7_4$ which is given as the plat closure of the braid word $\sigma_2^3 \sigma_1^{-1} \sigma_2^3$. Since the knot has determinant $15$, it is $3$-colorable. Since it is a $2$-bridge knot, 
the  $3$-fold branched cover of $S^3$ 
along the knot  can be constructed with three $0$ handles, two $1$-handles, two $2$-handles and three $3$-handles. The union of the $0$ and  $1$-handles is a $3$-ball as is the other side of the decomposition. So the covering space is also $S^3$. We construct, quite explicitly, an immersed folding. Yet, there is a simple closed curve of double points for this folding. We demonstrate why this knot does not bound an embedded folding. The result {\it does not} say that the $3$-sphere cannot be embedded, but instead says that an embedding does not project canonically upon the standard sphere as a covering. 

In the case of the $3$-fold branched covering of $S^4$ branched along the $2$-twist-spun trefoil, we have a folded immersion in $S^4\times D^2$, but not a folded embedding. This example will be presented elsewhere. 

Our main results are interesting in the light of Theorems of Alexander \cite{Al1920}, Hilden \cite{Hi1976} and Montesinos \cite{Mo1976} in dimension $3$ and Iori and Piergallini \cite{IP2002} in dimension $4$. 

\begin{theorem}[J.W. Alexander \cite{Al1920}]
For any closed oriented and connected $m$-manifold $M^m$, there exists a simple branched covering $f : M^m \to S^m$ for some degree.  
\end{theorem}

\begin{theorem}[H. M. Hilden \cite{Hi1976}, J. M. Montesinos \cite{Mo1976}] 
For any closed oriented and connected $3$-manifold $M^3$,  there exists a $3$-fold simple branched covering $f : M^3\to S^3$ such that  
 the branch set $L$ is a link (or a knot). 
\end{theorem} 

The following is a conjecture due to Montesinos.  

\begin{conjecture}
For any closed oriented and connected $4$-manifold  $M^4$,  there exists a $4$-fold simple branched covering
$f : M^4\to S^4$ such that $L$ is an embedded surface in $S^4$.    
\end{conjecture} 

Some partial answers to this conjecture are known as follows.

\begin{theorem}[R. Piergallini \cite{Pi1995}] 
For any closed oriented and connected $4$-manifold  $M^4$,  there exists a $4$-fold simple branched covering
$f : M^4\to S^4$ such that $L$ is an immersed surface in $S^4$.     
\end{theorem} 

\begin{theorem}[M. Iori and R. Piergallini \cite{IP2002}] 
For any closed oriented and connected $4$-manifold  $M^4$,  there exists a $5$-fold simple branched covering
$f : M^4\to S^4$ such that $L$ is an embedded surface in $S^4$.  
\end{theorem}

It is known that every closed $3$-manifold can be embedded in $\R^5$, and that the obstructions to embedding a closed orientable $4$-manifold in $\R^6$ are given by the Pontryagin class and the signature of the $4$-manifold.  So in the $4$-dimensional case, the surface of self-intersections will be related to one of these characteristic classes.

Complex projective space $\pm \Co P^2$ can be obtained as a $2$-fold branch cover of $S^4$ branched along the standardly embedded projective planes of normal Euler class $\pm 2$ \cite{Massey}. However, our results for $2$-fold branched covers are dependent upon the existence of a Seifert manifold for $2$-fold branched covers. The normal Euler classes prohibit the existence of such a Seifert solid. It is possible that an alternative chart movie description can handle this case, but as of this writing we do not have an alternative method for this. 

Our principal technique is to generalize the notion of a chart  to a chart surface (or curtain) in dimension $3$. In dimension $4$, we develop a $3$-dimensional analogue called an interwoven solid. We recall from \cite{Kam2002} that a {\it chart} is a labeled finite graph in the plane that has three types of vertices: $1$-valent black vertices, $4$-valent crossings, and $6$-valent white vertices. The labels upon the edges incident at crossings and $6$-valent vertices are required to satisfy additional conditions that we discuss below (Section~\ref{chartdef}).

\begin{wrapfigure}[14]{l}{4.2in}\includegraphics[width=4in]{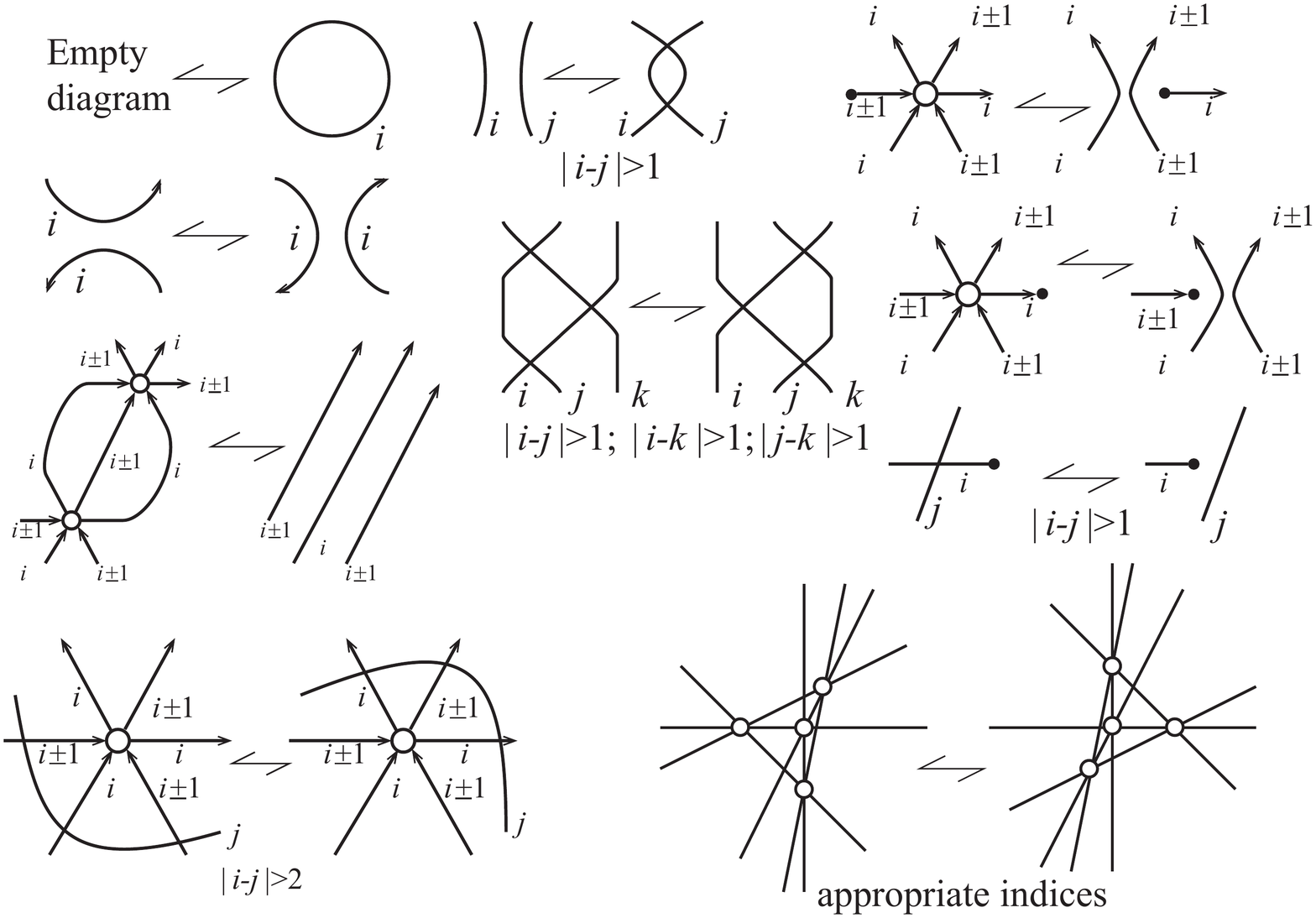}\end{wrapfigure}
A {\it curtain} is an immersed  labeled compact surface with boundary in $3$ space. The boundary is embedded as a knot or link which may pass through the sheets of the curtains. The singularities of the curtains are non-generic intersections among three sheets and transverse intersections between two sheets. The curtain can be put into general position with respect to a height function defined on $3$-space. In this case, the critical points are quantified as changes to charts. These changes are as follows:

\begin{enumerate}

\item introducing or removing a single labeled edge that joins two black vertices but does not intersect another edge of the chart;

\item an application of any of the {\it chart moves} that are depicted.

\end{enumerate}

A knotted or linked surface in $4$-dimensional space can be given via a movie description, for example. In such a movie a sequence of knot and link diagrams is given in such a way that successive cross-sections differ by a critical point (birth or death of a simple closed curve or a saddle point) or by one of the Reidemeister moves that are indicated.
\begin{wrapfigure}[10]{l}{3in}
\includegraphics[scale=.5]{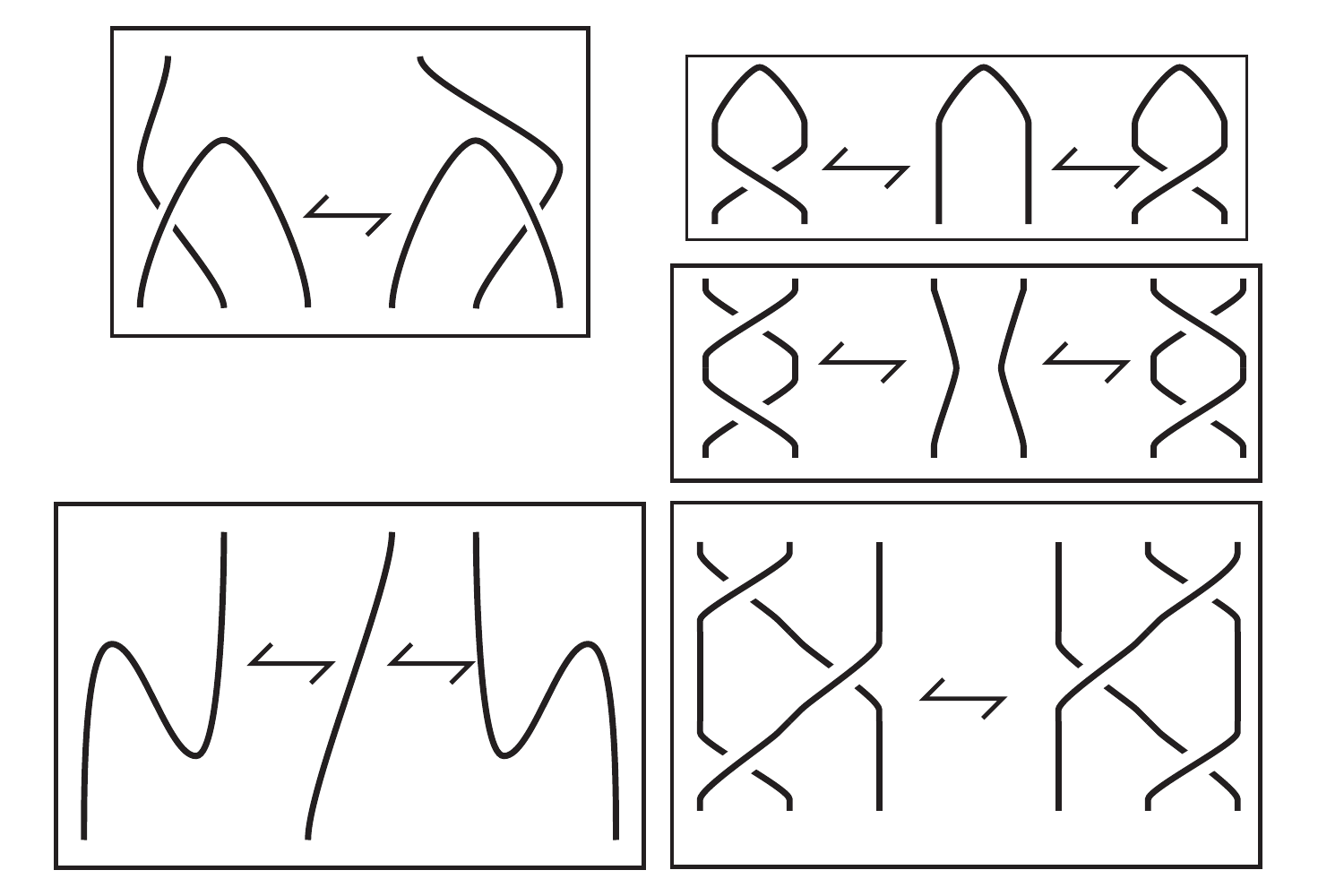} \vspace{-.5in}
\end{wrapfigure}
In such a movie, we explicitly choose a height function for each of the stills and keep this height function consistent throughout the movie. Alternatively, an oriented surface knot or link can be put into surface braid form (via an oriented chart), and the movie description is determined by a critical point analysis of the corresponding braid chart. We will have cause to use both descriptions.

In order to produce a folded immersion of the branched covering of $S^4$ branched along the surface link, we need a tertiary chart-like structure: that of an interwoven solid.

Interwoven solids are to curtains as curtains are to charts. Thus an {\it interwoven solid} is an immersed $3$-dimensional manifold in $4$-space whose boundary is a knotted or linked surface. The boundary may pass through the solid, and the immersion is not in general position. An interwoven solid can be described as a sequence of curtains (with nodes) where successive curtains differ by generic critical points or specific curtain moves. Here we will not list all the possible curtain moves, but leave that taxonomy for our future work in this direction. We give three specific examples of interwoven solids. First, the Seifert solid that a knotted or linked surface bounds in $4$-space is an interwoven solid. Second, consider a genus $n$ surface that is  embedded in $3$-space. An interwoven solid is the bounded $3$-manifold that such a surface bounds. Third, our last example consists of the spun trefoil presented in movie form. The interwoven solid that yields its $3$-fold branched cover is indicated as a sequence of curtains. Successive curtains in this case differ from each other by replacing a chart move and its inverse with a product of charts or {\it vice versa}.

So for each still in the movie (or braid movie) of the knotted surface, we construct a curtain. In the case of the $2$-fold branched cover, the curtain is  a Seifert surface. Each critical event for the knot movie induces a critical event between the curtains. Each will be described explicitly. At the end of the movie, it may be necessary to unravel the curtains. In our examples for the $3$-fold branch covers of classical knots, we need to use chart moves to achieve the unraveling. 

Here we will introduce some of the curtain movie moves that are necessary to achieve the unraveling of the 
interwoven solid. However, a full list is not necessary to articulate since we can, at the end of the movie, isolate the curtains into an immersed closed surface that is contained in a ball neighborhood in a $3$-dimensional slice of $4$-space. This surface can be eliminated since it is null-homologous in $3$-space. 

Here is an outline of the paper. We begin with a basic review of the permutation and braid groups. Then we describe simple branched covers of $2$-dimensional surfaces. We extend the ideas of the second author to give a combinatorial description of permutation and braid charts. In particular, we demonstrate that vertices in charts correspond to specific critical events for surfaces in $3$ and $4$-dimensional space. We then use our ideas to develop embeddings of $2$-fold branched covers of the $3$ and $4$-sphere branched over codimension $2$ embedded submanifolds. Many examples are given with quite a bit of detail rendered. Three examples of embedded and immersed  foldings of the $3$-fold branched covers of the $3$ and $4$ sphere are given. The paper closes with a proof of Theorem~\ref{mainimmersed} in case $k=1$.

\subsection*{Acknowledgements.} The genesis of this paper was  a conversation that JSC had with John Etnyre at Georgia Tech a few years ago. During the sixth East Asia Knot School in Hiroshima (January 2011), JSC and SK improved upon the original construction. At that time JSC had support from NSF grant \#0603926. JSCs visit to Hiroshima was generously supported by grant JSPS KAKENHI  number 19204002 issued to Makota Matsumoto.   This paper began when JSC was visiting Kyungpook National University. His visit was  supported by the Ministry of Education Science and Technology (MEST) and the Korean Federation of Science and Technology Societies (KOFST). SK is being supported by  JSPS KAKENHI Grant  number  21340015 and 23654027. As part of that grant, he was able to visit Kyungpook, continue this conversation with JSC, and the authors developed the final steps in the process.  We also would like to thank Makoto Sakuma and Daniel Silver for helpful bibliographical remarks.

Let us proceed.

\section{Some basic definitions}

Throughout this work, we will be using the permutation and braid groups. The permutation group $\Sigma_n$ is the set of bijective maps on $\{1,2,\ldots, n\}$. It has a presentation:
\begin{eqnarray*}\Sigma_n= \langle \tau_1, \tau_2, \ldots, \tau_{n-1} &: & \tau_i \tau_j = \tau_j \tau_i \ {\mbox{\rm if }}\ |i-j|>1; \\ && \tau_i \tau_{i+1} \tau_i = \tau_{i+1} \tau_i \tau_{i+1} \ {\mbox{\rm if }} \ i=1,\ldots n-2;  \\  && \tau_i^2 =1  \ {\mbox{\rm if }} \ i=1,\ldots, n-1 \rangle. \end{eqnarray*}
The permutation group is a quotient of the braid group $B_n$ which has the presentation:
\begin{eqnarray*} B_n = \langle \sigma_1, \sigma_2, \ldots, \sigma_{n-1} &: & \sigma_i \sigma_j = \sigma_j \sigma_i \ {\mbox{\rm if }}\ |i-j|>1; \\ && \sigma_i \sigma_{i+1} \sigma_i = \sigma_{i+1} \sigma_i \sigma_{i+1} \ {\mbox{\rm if }} \ i=1,\ldots n-2 \rangle.  \end{eqnarray*}
We will be considering knotted and linked subsets of the {\it  $(k+2)$-dimensional sphere} for $k=-1,0,1,2$ which is defined to be the set of unit vectors in $\R^{k+3}$. Specifically, $$S^{k+2}= \{ (x_1,x_2, \ldots , x_{k+3}) \in \R^{k+3}: \sum_{j=1}^{k+3} x_j^2 = 1 \}.$$ This is the boundary of the unit $(k+3)$-disk
$$D^{k+3}= \{ (x_1,x_2, \ldots , x_{k+3}) \in \R^{k+3}: \sum_{j=1}^{k+3} x_j^2 \le 1 \}.$$ 
We will frequently consider a fixed embedding of a $(k+2)$-dimensional disk in the sphere $S^{k+2}$. Thus many constructions occur in the disk and extend trivially outward to the sphere. We include $\R^{k+3}$ in $\R^{k+4}$ in a standard fashion as $\{(x_1,x_2, \ldots , x_{k+3}, 0): x_j\in \R \}$.
Furthermore, It will often be convenient to decompose the disk as a cartesian product of lower dimensional disks. For example, $D^2 \simeq D^1 \times D^1 \simeq I \times I$, where $I$ is the unit interval $I=[0,1]$. 

A {\it simple branched cover of $S^{k+2}$ of degree $n$} is a compact orientable manifold $M^{k+2}$ together with a surjective map $f:M^{k+2} \rightarrow S^{k+2}$ such that each point in $S^{k+2}$ is of one of two types.

\begin{enumerate}
\item  A {\it regular point} is a point that has a $(k+2)$-disk neighborhood $N$ such that $f^{-1}(N)$ is homeomorphic to the disjoint union of exactly $n$ copies of $N$.
\item A {\it (simple) branch point} is a point that has a neighborhood $N$ such that $f^{-1}(N)$ consists of $n-1$ disk neighborhoods $N_1, \ldots, N_{n-1}$, the map 
$f|_{N_j}$ 
is a homeomorphism for $j=1, \ldots, n-2$, and $f|_{N_{n-1}} : N_{n-1} \rightarrow N$ is a two-to-one branched covering map. In this case, there is a coordinate system $(x_1, \ldots , x_{k+2})$ of this component such that the restriction of $f$ is given by $(x_1,x_2, x_3, \ldots, x_{k+2}) \mapsto (x_1^2-x_2^2, 2x_1x_2, x_3 ,\ldots, x_{r+2})$.
\end{enumerate}

A branched cover of degree $n$ is also called an {\it $n$-fold branched cover}.
Throughout, we will only deal with simple branched coverings, and so we will speak colloquially of branched coverings. We sometimes work in the PL-category.

\section{$2$-dimensional simple branched coverings}

Let 
$f: M^2 \to S^2$ be an $n$-fold simple branched cover with branch set $L$, and let 
$\underline f : M^2 \setminus f^{-1}(L) \to S^2 \setminus L$ be the associated covering map; that is $\underline{f}$ is the restriction of $f$ to the complement of the branch set.  

Take a base point $\ast$ of $S^2 \setminus L$ to consider the fundamental group $\pi_1(S^2 \setminus L, \ast)$. 
The preimage $f^{-1}(\ast)$ of the base point $\ast$ consists of $n$ points of $M^2$.  
Then we have a  {\it monodromy } $ \rho: \pi_1(S^2 \setminus L, \ast) \to \Sigma_n$, where the symmetric group $\Sigma_n$  on letters $\{1, 2, \dots, n\}$  is identified with the symmetric group on $f^{-1}(\ast)$.   (A monodromy $\rho$ depends on the identification between  $\{1, 2, \dots, n\}$ and $f^{-1}(\ast)$.)   The covering map $\underline f $ is determined by the monodromy.  

By the Riemann-Hurwitz formula,  $L$ consists of an even number of points.  The idea of the monodromy is illustrated in the central figure below. The concept of a (permutation) chart was summarized in the introduction and is discussed in detail below.

\begin{center}
\includegraphics[width=4in]{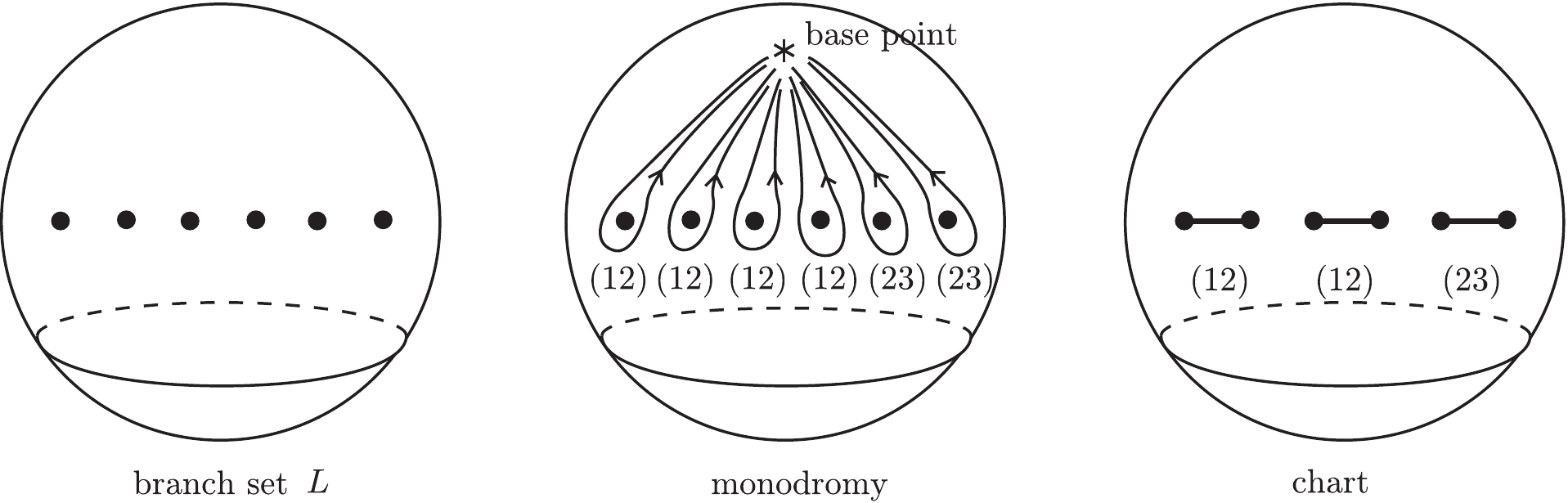}
\end{center}

\begin{center}
\includegraphics[width=4in]{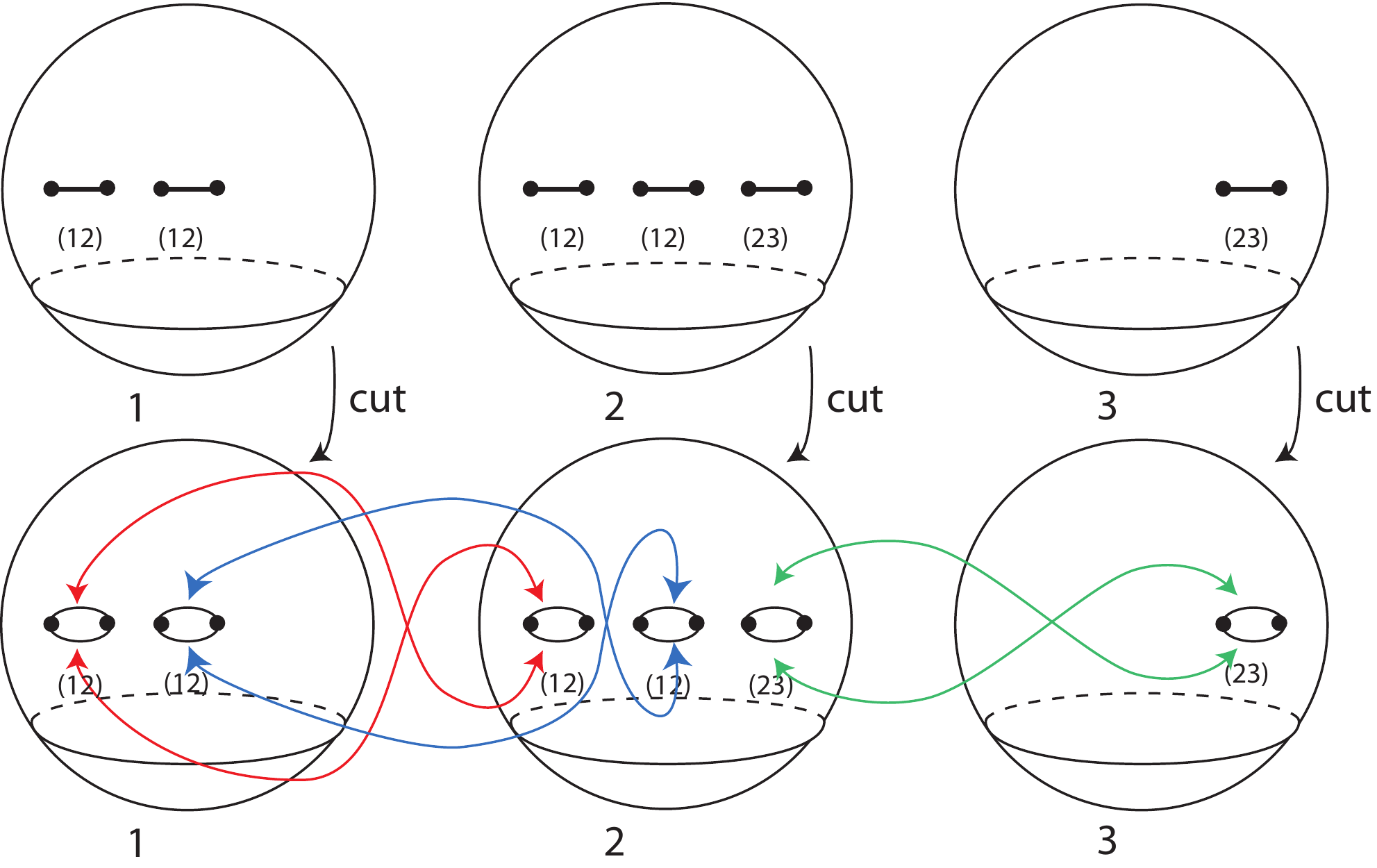}
\end{center}

\begin{wrapfigure}[10]{l}{1.4in}\vspace{-.2in}
\includegraphics[width=1.3in]{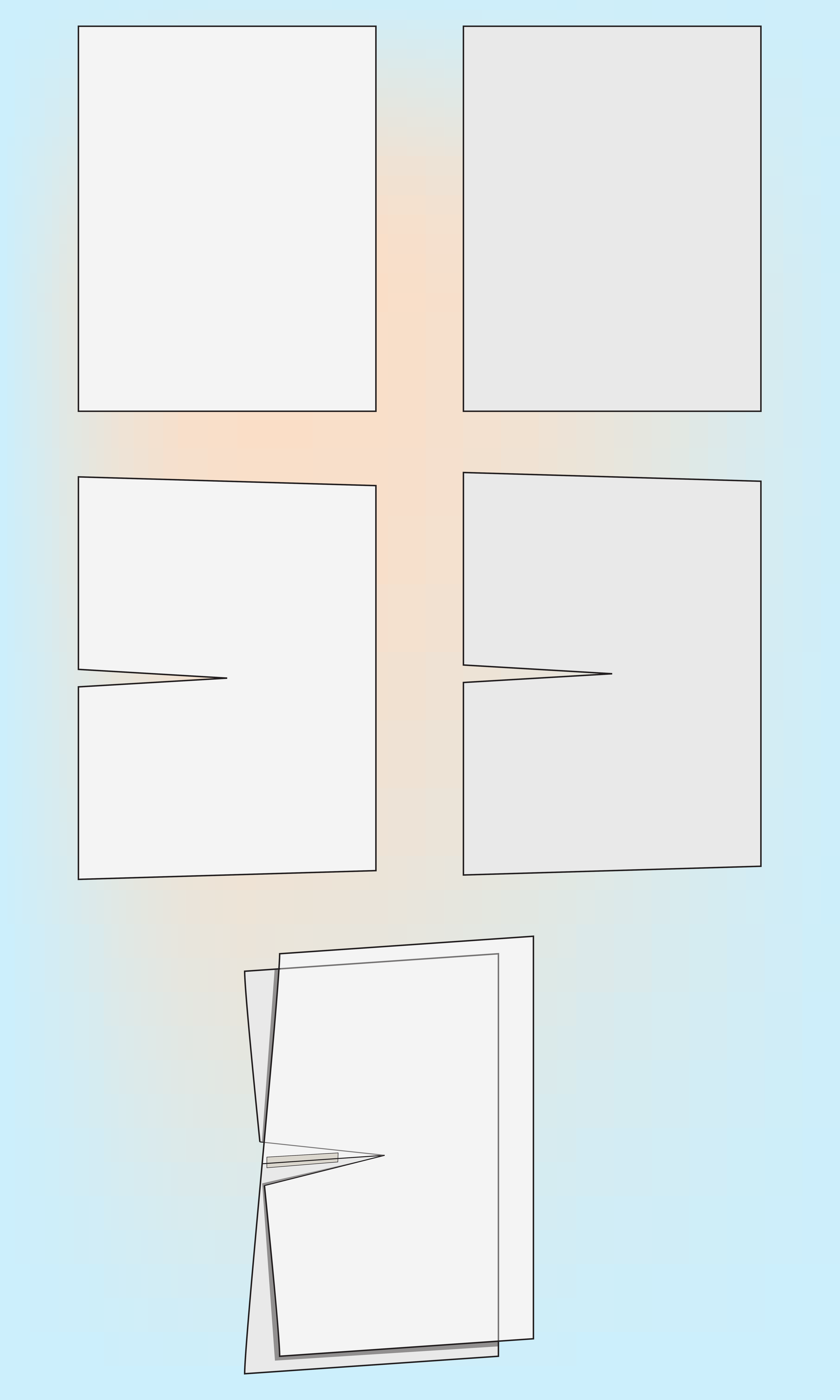}
\end{wrapfigure} When a monodromy is described by a chart, it is easy to 
construct $M^2$. We explain it by using an example.   
Let $\Gamma$ be the chart depicted on the right of the illustration above. 
Consider three copies of $S^2$ labeled by $1$, $2$, and $3$, say $S^2_1$, $S^2_2$ and $S^2_3$, respectively.    
On the copy $S^2_1$, draw the edges with label $(12)$ of $\Gamma$, 
on the copy $S^2_2$, draw the edges with label $(12)$ of $\Gamma$ and those with label $(23)$, 
and on the copy $S^2_3$, draw the edges with label $(23)$.  
Cut the three $2$-spheres along these edges, to obtain three compact surfaces, 
say $M_1$, $M_2$ and $M_3$, as indicated here.

The surface $M^2$ is obtained from the union $M_1 \cup M_2 \cup M_3$  by identifying the boundary in a way as follows:  Let $e$ be an edge 
with label $(12)$ on $S^2_1$, and let $e_+$ and $e_-$ be the copies of $e$ in $\partial M_1$.  
Let $e'$ be the corresponding edge on $S^2_2$, and let $e'_+$ and $e'_-$ be the corresponding copies in $\partial M_2$.  Then we identify $e_+$ with $e'_-$, and identify $e_-$ with $e'_+$, respectively.  All boundary edges of $M_1 \cup M_2 \cup M_3$ are identified in this fashion, and we have a closed surface.  This is the desired $M^2$.  

\begin{wrapfigure}[4]{r}{1.7in} \vspace{-.4in}\includegraphics[width=1.7in]{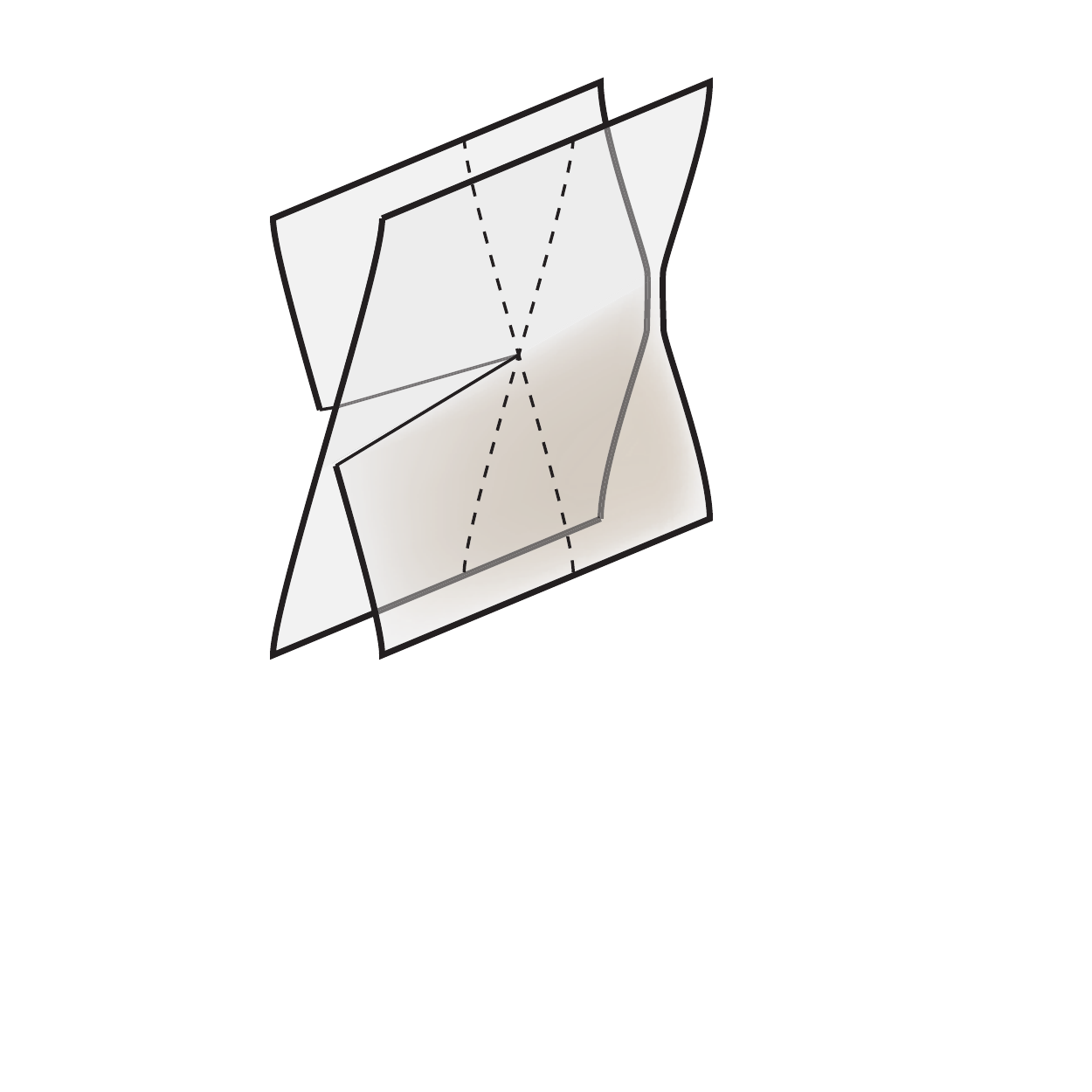}\end{wrapfigure} In neighborhood of the branch points we can cut and paste (tape) as indicated with two sheets of paper here. The local picture of a branch point as a broken surface diagram is indicated.

\subsection{Permutation and braid charts}
\label{chartdef}

A {\it (permutation) chart of degree $n$}, (or  an {\it $\Sigma_n$-chart}), is a labeled finite graph embedded in the $2$-disk $D^2 \simeq I \times I$ which has three types of vertices. The labels are written on the edges are taken from the set $\{1,2, \ldots , n-1 \}$ --- these will correspond to the generators $\tau_1, \ldots, \tau_{n-1}$ of the permutation group $\Sigma_n$. The vertices are of the following type:
\begin{enumerate}
\item a {\it black vertex} is a mono-valent whose incident edge may have any label taken from  $\{1,2, \ldots , n-1 \}$;
\item a {\it crossing} is $4$-valent, and the labels on the incident edges are given in cyclic order $i,j,i,j$ where $|i-j|>1$;
\item a {\it white vertex} is $6$-valent, and its incident edges have labels in cyclic order given by $i,i+1,i,i+1,i,i+1$ for some $i =1,2,\ldots, n-2$. 
\end{enumerate}
Necessarily, a chart has an even number of black vertices. 

Let $\Gamma \subset D^2$ denote a chart with $m$ black vertices $\{b_1,\ldots, b_m\}$. We identify $D^2$ with $[0,1] \times [0,1]$ and consider the projection $p_1:[0,1]\times [0,1] \rightarrow [0,1]$ onto the first factor. We will assume that $\Gamma$ is in {\it  general position with respect to $p_1$}. That is, the critical points of the arcs are all non-degenerate ($C^2$-approximated by quadratic functions), and each critical points  or vertex projects to a different time value. For convenience, we describe vertices as critical points of the chart.  
We may order the black vertices from left-to-right with respect to $p_1$. In this way the image of $b_j$ under the projection $p_1$ onto the first factor is $t_j$ and $0< t_1 < t_2 < \cdots < t_m < 1.$

By the correspondence $i  \leftrightarrow  \tau_i = (i,  i+1) \in \Sigma_n$, the labels of a chart are assumed to be 
transpositions in $\Sigma_n$. For a chart $\Gamma$, we consider a monodromy 
$$\rho_\Gamma: \pi_1(D^2 \setminus L) \to \Sigma_n,  \quad [\ell] \mapsto [  \textrm{intersection word of $\ell$ w.r.t. $\Gamma$}   ],  $$ 
where $L$ $(= L_\Gamma)$ is the set of black vertices $\{b_1, \ldots, b_m\}.$  An intersection word is a sequence of elements of  $\{1, \dots, n-1\}$, which is regarded  as an element of $\Sigma_n$ by the correspondence $i  \leftrightarrow  \tau_i = (i ~  i+1) \in \Sigma_n$.

\begin{wrapfigure}[12]{r}{3in}
\includegraphics[scale=0.12]{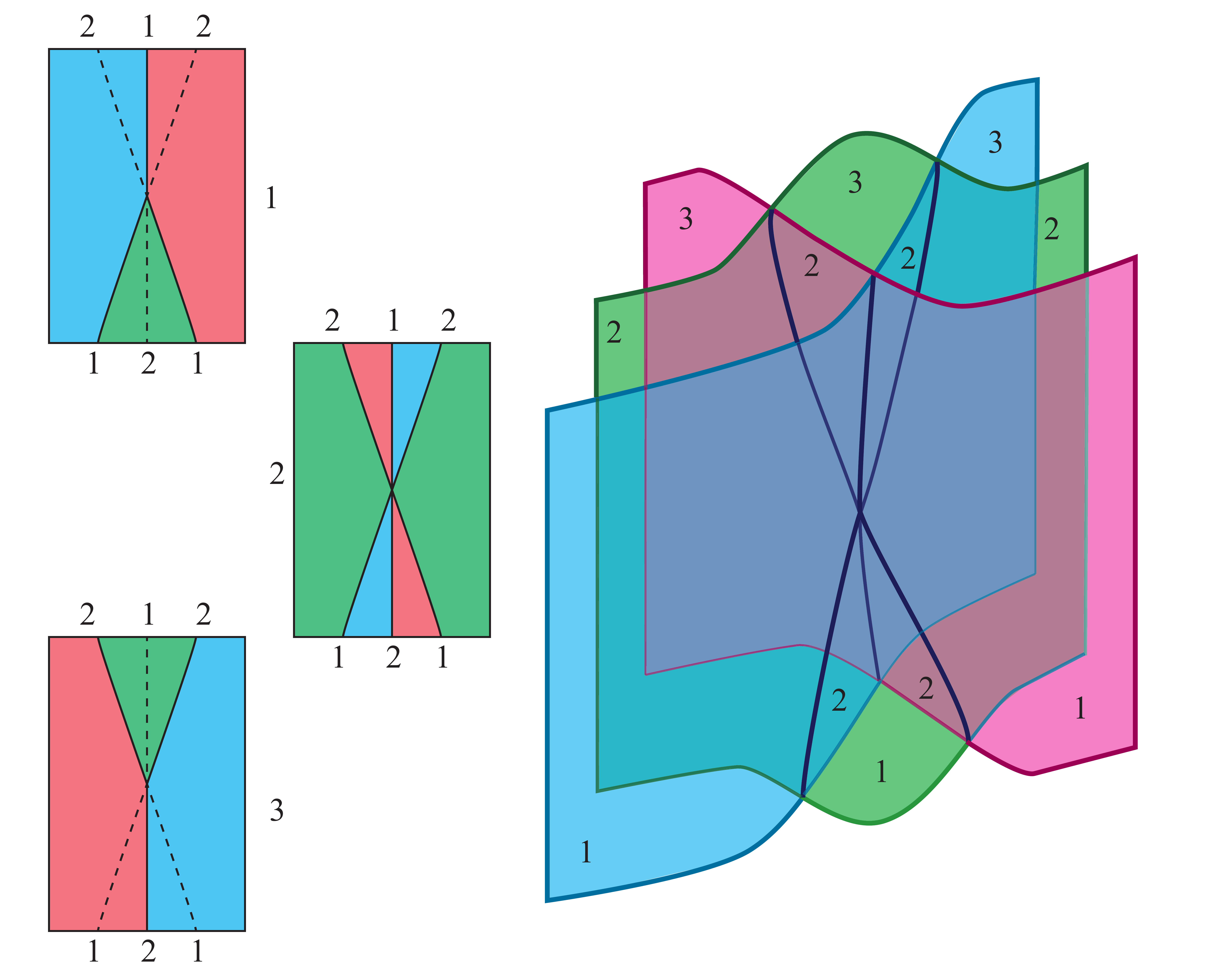}
\end{wrapfigure}  
If an edge terminates at a black vertex, then this becomes a simple branch point as above. If an edge terminates at a crossing, the $i\/$th and $(i+1)\/$st sheets are cut and re-glued while the $j\/$th and the $(j+1)\/$st sheets are cut and re-glued. 
The identifications occur in disparate sheets and the cutting can proceed. In the case that an edge terminates at a white vertex the reassembly is indicated.
The resulting configuration resembles the intersection of the three coordinate planes in space. The resulting surface $M^2(\Gamma)$  with boundary is called the {\it $n$-fold irregular  simple branched cover induced from the chart\/  $\Gamma$}. The surface $M^2(\Gamma)$ extends to a simple branched cover of the $2$-sphere, by trivially extending the cover over the complement of the disk in $S^2$.

The classification of simple branched coverings was studied by 
J. L{\" u}roth \cite{Lur1871}, A. Clebsch \cite{Cl1973} and  A. Hurwitz \cite{Hur1891}. 
The classification theorem is stated as follows.  

\begin{theorem}\label{thm:classification1}
Let $f : M^2 \to S^2$ and $f' : {M^2}' \to S^2$ be   $n$-fold simple branched coverings with branch sets  $L$ and $L'$, respectively.  We assume that $M^2$ and ${M^2}' $ are connected.  Then 
$ f $ and $ f'$ are equivalent if and only if $\# L = \# L'$.  
\end{theorem}

Hurwitz  \cite{Hur1891} studied branched coverings by using of a system of monodromies of meridian elements of the branch set, called a 
{\it Hurwitz system}, and studied when two systems present the same (up to equivalence) branched coverings.  

A Hurwitz system depends on a system of generating set of $\pi_1(S^2 \setminus L, \ast)$.  
For the generating system depicted in the middle of the first illustration of this section, the 
Hurwitz system is  
$$\alpha = ((12), (12), (12), (12), (23), (23)).$$   
Besides a choice of a generating system, a Hurwitz system  depends  on the identification of 
$\{1, 2, \dots, n\}$ and the fiber $f^{-1}(\ast)$.  

Two Hurwitz systems present the same (up to equivalence) monodromy if and only if they are 
related by a finite sequence of  {\it Hurwitz moves}  and {\it conjugations}.  
The {\it Hurwitz moves} are 
$$(a_1, \dots, a_k, a_{k+1}, \dots, a_m) \mapsto (a_1, \dots, a_{k+1}, a_{k+1}^{-1} a_k a_{k+1}, \dots, a_m)$$ 
for $k=1, \dots, n-1$ and their inverse moves.   {\it Conjugations} are 
$$(a_1, \dots,  a_m) \mapsto (g^{-1} a_1 g, \dots, g^{-1} a_m g)$$
for $g \in \Sigma_n$.   When two Hurwitz systems are related by a finite sequence of  {Hurwitz moves  and conjugations, we say that they are {\it HC-equivalent}.  ($H$ and $C$ stand for Hurwitz and conjugation.) 

Due to Hurwitz  \cite{Hur1891}, the classification theorem is stated as follows.  

\begin{theorem}
Let $f : M^2 \to S^2$ be an $n$-fold simple branched covering.  Assume that $M^2$ is connected.  
Any Hurwitz system of $f$ is HC-equivalent to 
$$((12), \dots, (12), (13), (13), (14), (14), \dots, (1,n), (1,n)).$$ 
\end{theorem}

\begin{theorem} 
Let  $f : M^2 \to S^2$ be an $n$-fold simple branched covering, and  $\rho_f$ a monodromy of $f$.  
There exists a chart $\Gamma$ such that $\rho_\Gamma = \rho_f$.  
{\rm (We call $\Gamma$ a {\it chart description} of $f$ or $\rho_f$.)}
\end{theorem} 
{\sc Proof.} For a detailed proof, see \cite{Kam2002}. $\Box$

\subsection{Constructing a branched covering from a chart}

\begin{wrapfigure}[5]{r}{2in}
\vspace{-.35in}
\includegraphics[scale=0.35]{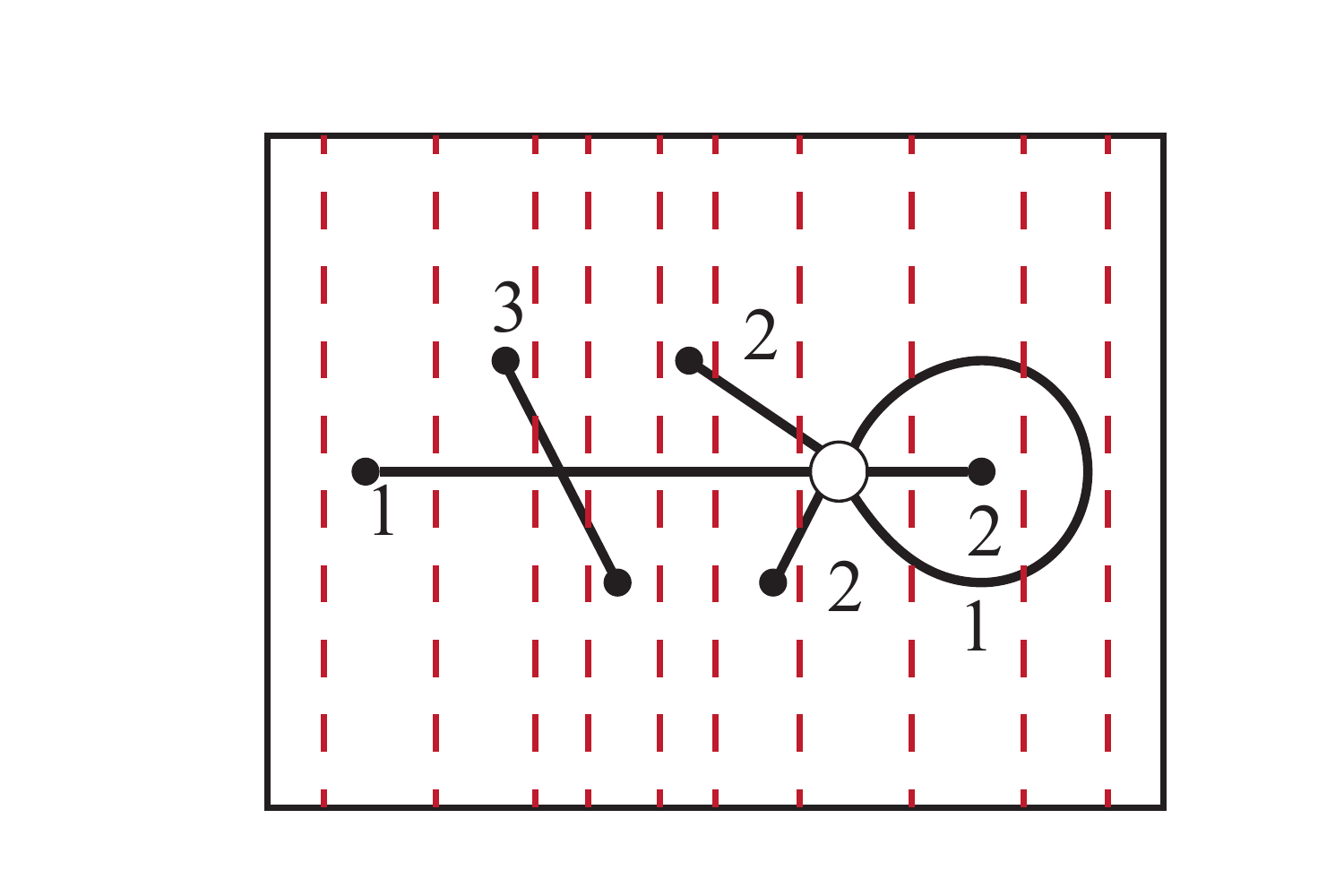}
\end{wrapfigure}
The surface $M^2(\Gamma)$ can be mapped into $3$-dimensional space in general position by using the chart to create a {\it permutation movie} that we turn to describe now. Recall, that the chart is in general position with respect to the projection $p_1$ onto the horizontal edge of the square. 
Each critical point of an edge is non-degenerate, and the critical points and vertices occur at distinct times. Let 
$t_0= 0<t_1 < t_2< \cdots < t_{s-1} <1 =t_s$ denote the critical values,  and consider the intersection of $\Gamma$ with the vertical line $L_i=p_1^{-1}((t_{i}+t_{i+1})/2).$ 

\begin{wrapfigure}[7]{l}{3in}\vspace{-.3in}
\includegraphics[scale=0.05]{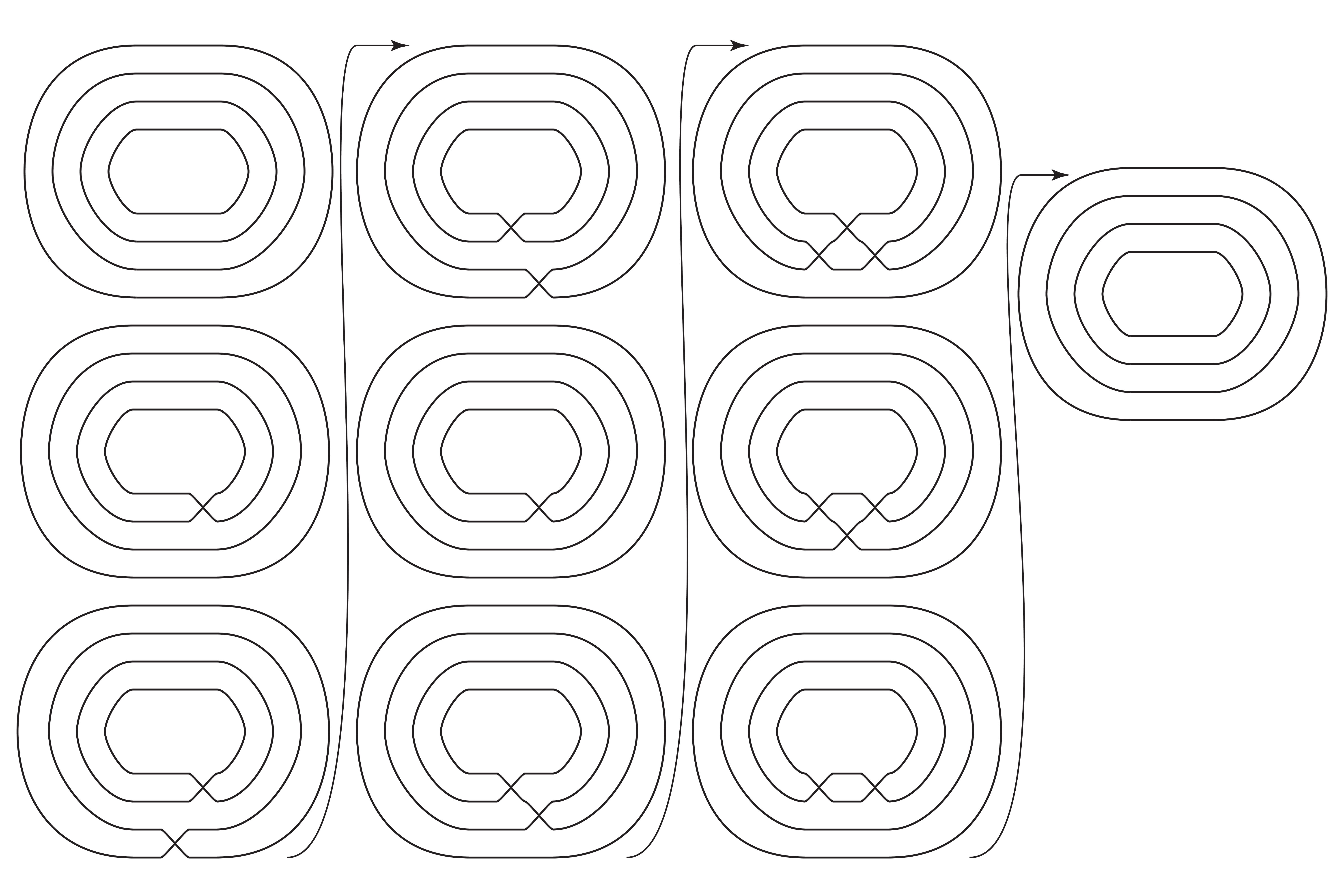}\end{wrapfigure} 
Now read the intersection sequence of $\Gamma$ with the lines $L_i$ for $i=1, \ldots , s$ in order. We read the sequence from top to bottom and write the sequence from left to right. For example, for the chart  above,the sequence reads: $\emptyset;$ $(1);$ $(3,1);$ $(1,3);$ $(1);$ $(2,1);$  $(2,1,2);$ $(1,2,1);$  $(1,1);$ $\emptyset.$ 
We rewrite this as a sequence of words in the permutation group $\Sigma_n$ (here $n=3$) as follows:
$1;$  $ \tau_1;$  $\tau_3\tau_1;$  $\tau_1\tau_3;$ $ \tau_1;$ $\tau_2\tau_1;$ $\tau_2\tau_1\tau_2;$  $\tau_1\tau_2
\tau_1;$ $\tau_1\tau_1;1.$

The critical events of  a permutation chart correspond to familiar pieces of surfaces that are mapped into $3$-space.  

\begin{wrapfigure}[5]{r}{2in}\vspace{-.4in}
\includegraphics[scale=0.12]{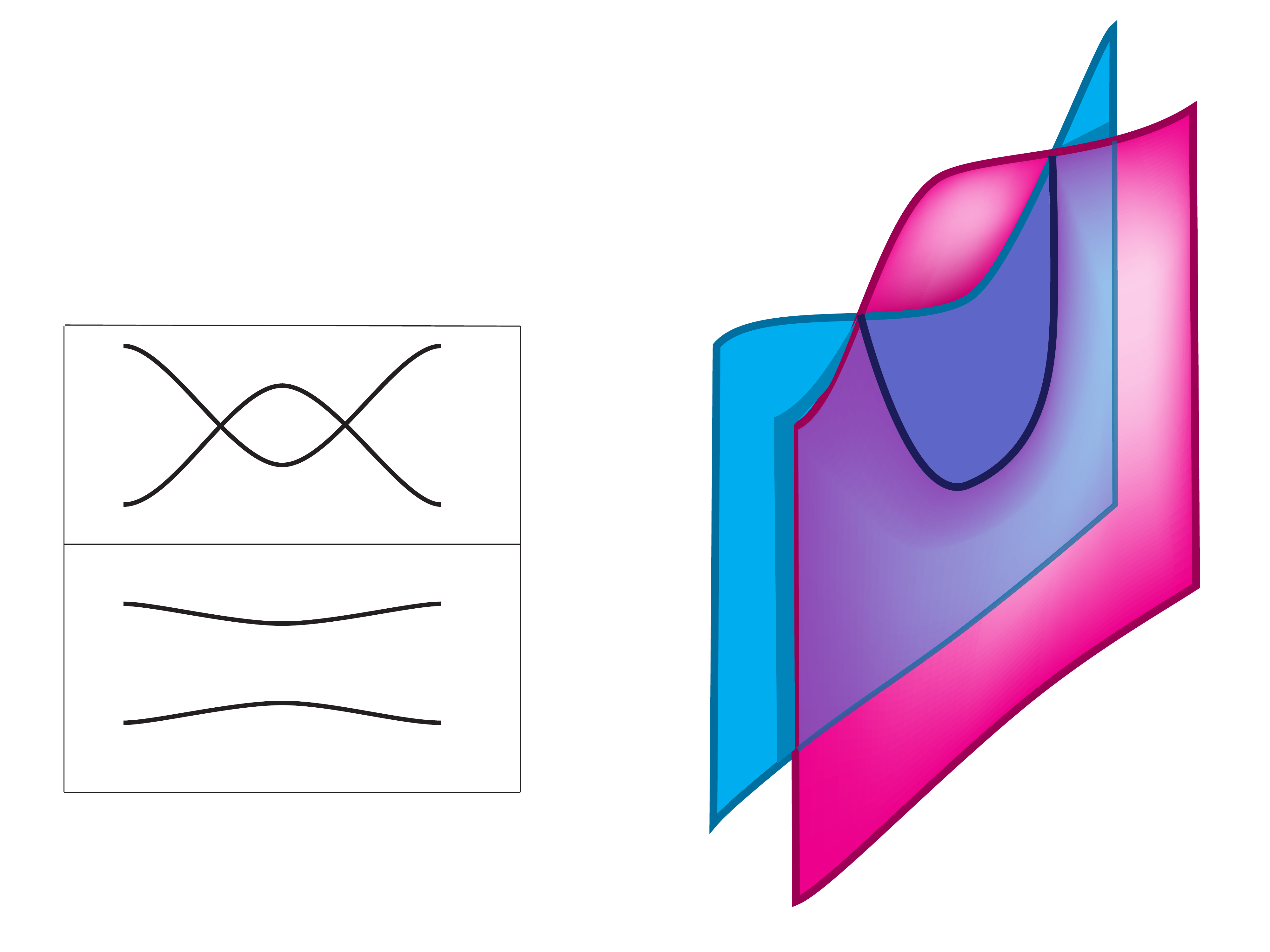}
\end{wrapfigure}
\noindent
{\bf 1. Type II moves.}
Critical events of the form $\subset^k_k$ correspond to the movie in which the creation of a pair of canceling generators $1 \Rightarrow \tau_k \tau_k$ occurs in a type II birth fashion. Critical events of the form ${}^k_k\supset$ correspond to the annihilation of the same pair of generators in a type II death fashion: $ \tau_k \tau_k \Rightarrow 1 $.  

\vspace{12pt}
\begin{wrapfigure}[4]{l}{2in}\vspace{-.7in}
\includegraphics[scale=0.3]{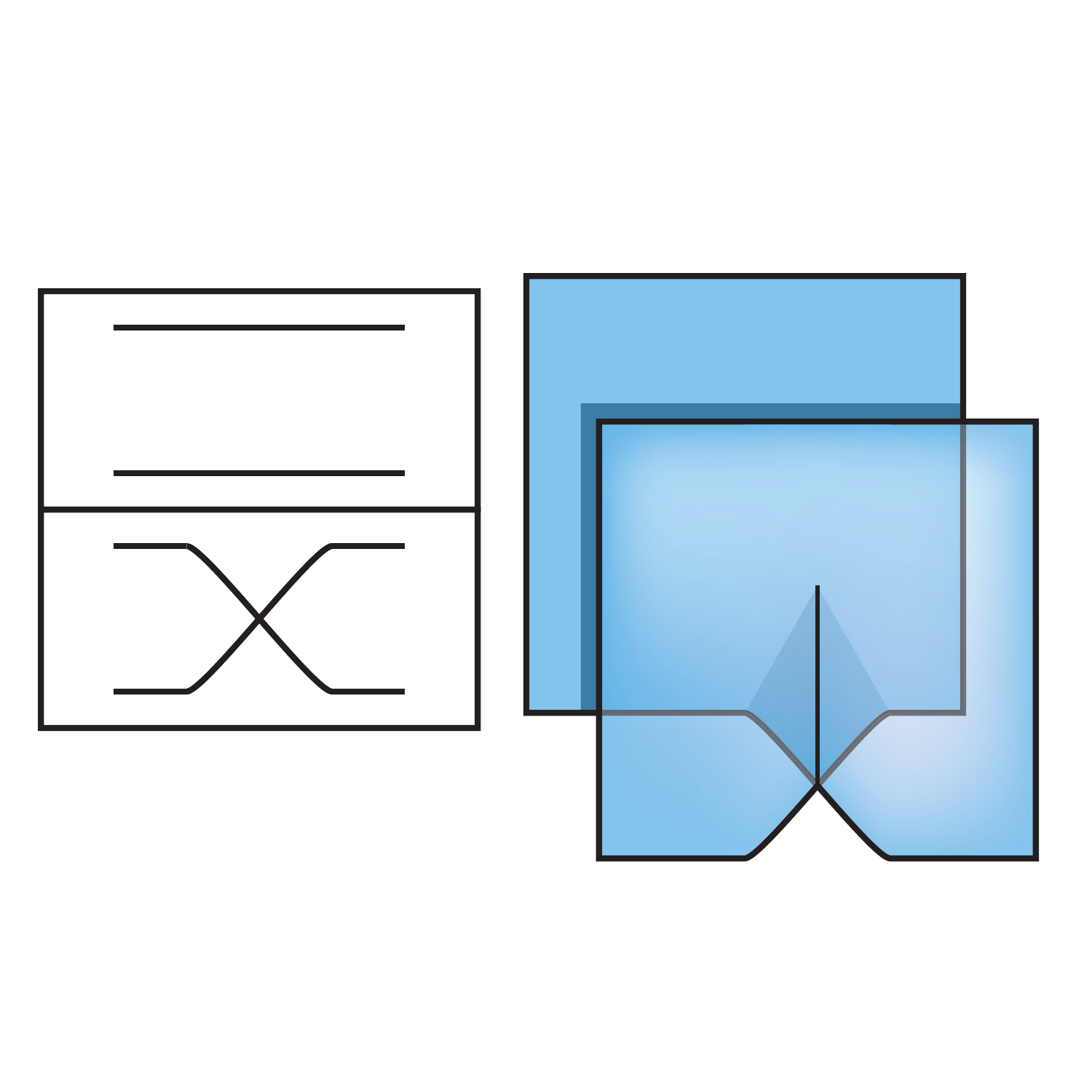}
\end{wrapfigure}
\noindent
 {\bf 2. Black vertices.} A black vertex corresponds to inserting or removing a generator  from a word. In the local picture of the surface, a branch point occurs at a saddle point in the surface.

\vspace{-.4in}

\noindent
{\bf  3. Crossing exchange.}
 A crossing indicates the interchange of distant crossings. 
 \rule{0in}{.75in}
 \rule{2in}{0in}\includegraphics[scale=0.14]{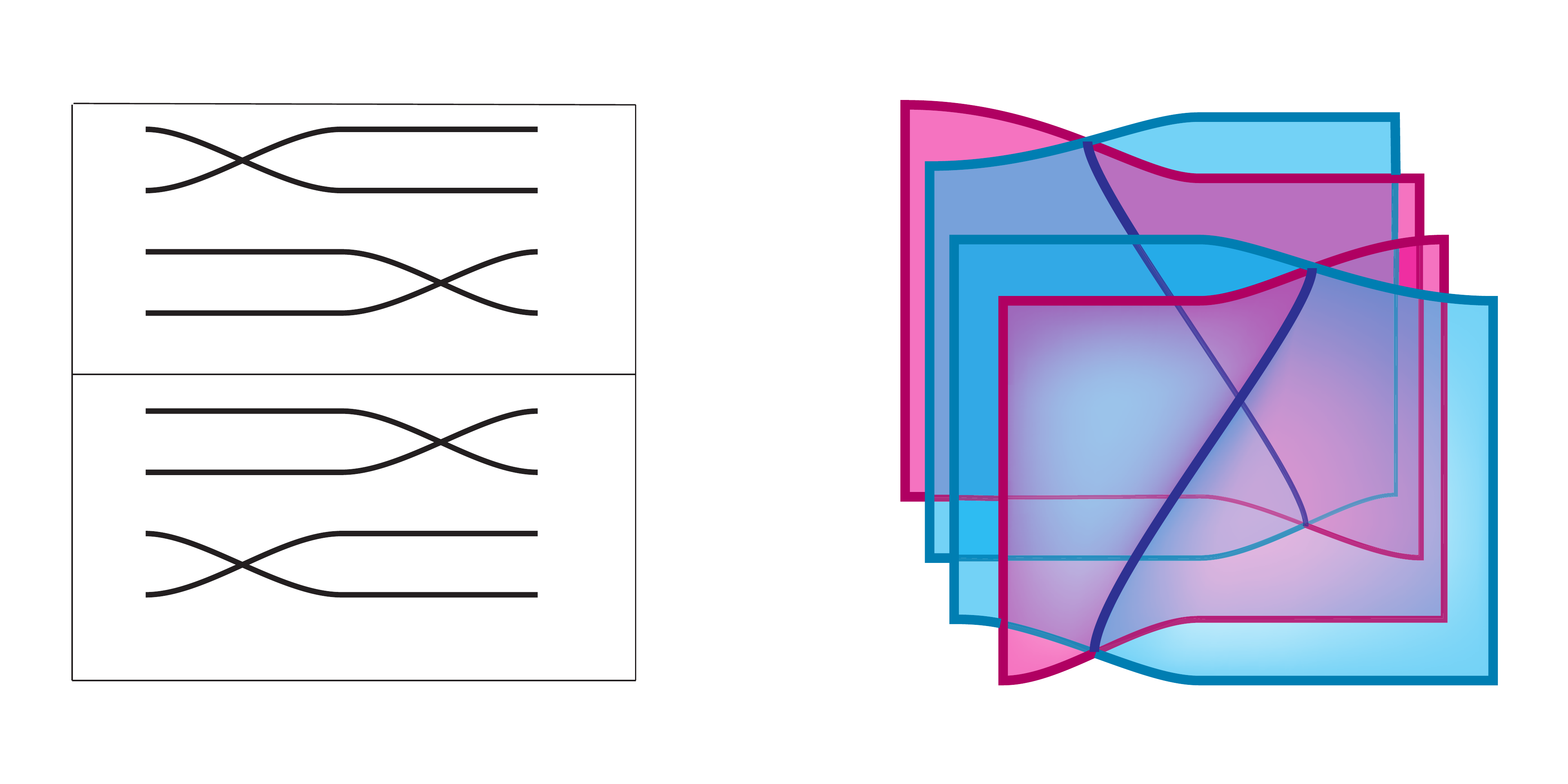}

\vspace{0.1in}

 \begin{wrapfigure}[4]{l}{2.2in}\vspace{-.5in}
\includegraphics[scale=0.07]{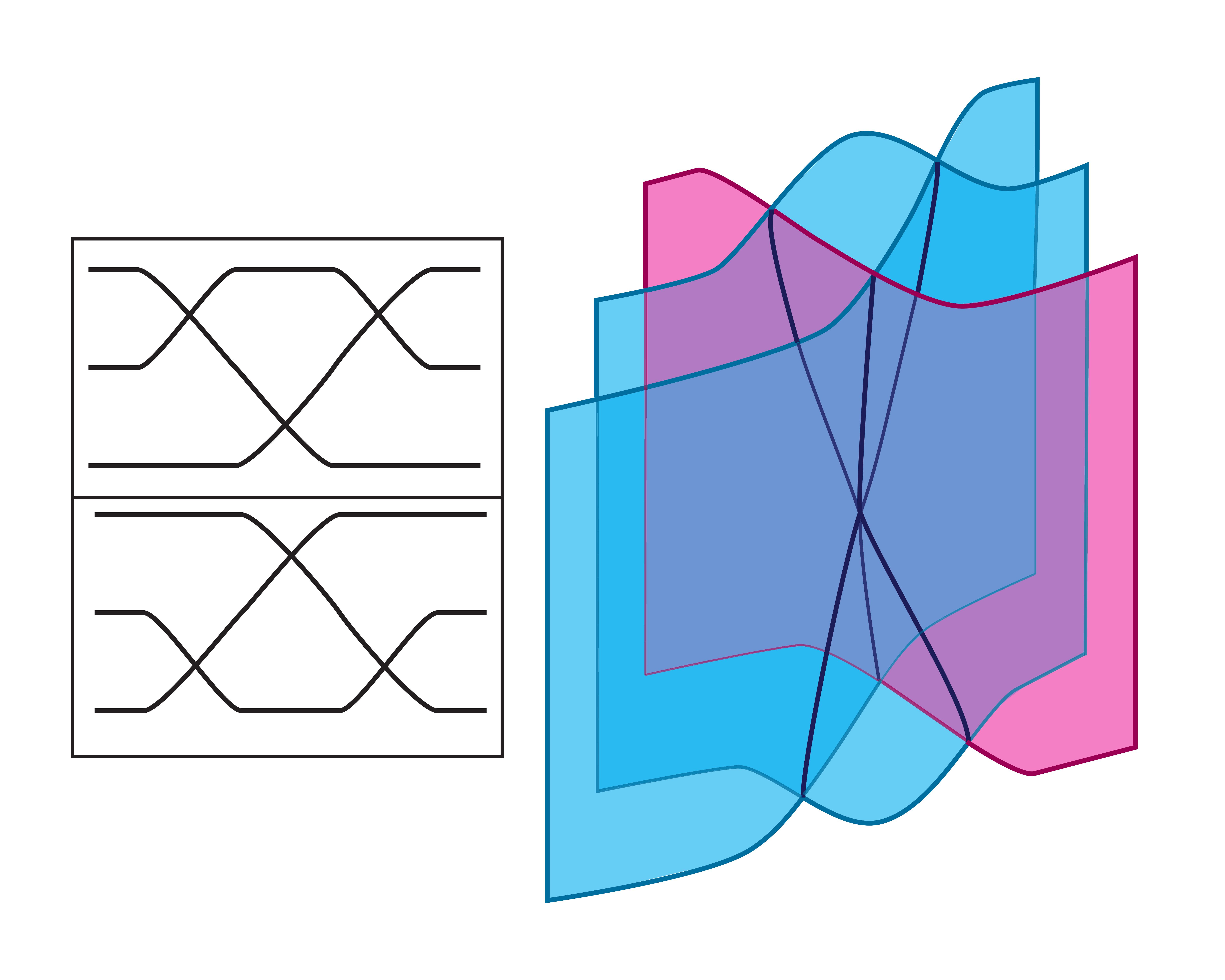}
\end{wrapfigure}
\noindent
{\bf 4. White vertices.}
A white vertex indicates a triple point of the surface induced by the relation $\tau_i \tau_{i+1}\tau_i = \tau_{i+1}\tau_i \tau_{i+1}.$

We have the following construction.

\begin{proposition} \label{perm} Let $\Gamma$ denote a permutation chart of degree $n$ in the disk $D^2$. Let $M^2(\Gamma)$ denote the $n$-fold irregular branched cover of $D^2$ induced by $\Gamma$. Then there is a general position map $\tilde{f}:M^2(\Gamma) \rightarrow D^2 \times I$ such that the projection onto $D^2$ induces the branched covering map. 

Let $F: \hat{M}^2(\Gamma)\rightarrow S^2$ denote the extension of the branch cover $M^2(\Gamma)$ over $D^2$ to the sphere $S^2$. Then there is an general position map $\tilde{F}: \hat{M}^2(\Gamma)\rightarrow S^2 \times [0,1]$ such that $p\circ \tilde{F} = F$ where $p$ is the projection onto the first factor. 
\end{proposition}

\subsection{Lifting a permutation movie to a braid movie}

\begin{wrapfigure}[9]{l}{1.2in}\vspace{-0.1in}
\includegraphics[scale=0.4]{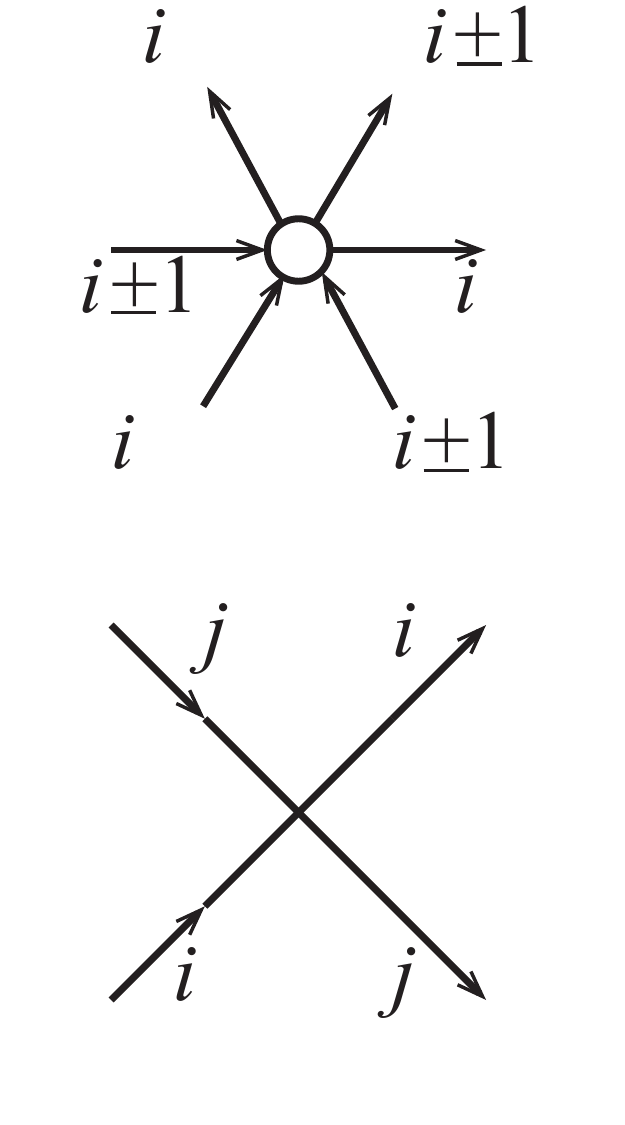}
\end{wrapfigure} A {\it braid chart of degree $n$} is a permutation chart in which orientations on the edges have been chosen so that  white vertices 
and crossings 
are of the form indicated 
in the figure. 
Thus there are  three    incoming and three outgoing edges that alternate $i$, $i\pm 1$, $i$, $i \pm 1$, $i$, $i \pm 1$ in cyclic order.  Similarly, at a crossing the edges with labels $i$ and the edges with labels $j$ are oriented consistently ($|i-j|>1$) as indicate in the figure. In this case, we say that there is a {\it flow though the white vertex} or {\it the crossing.} If a permutation chart can be consistently oriented to create a braid chart, then the permutation movie can be lifted to a braid movie. Specifically, right pointing edges with label $i$ correspond to braid generators $\sigma_i$, and left pointing edges with label $i$ correspond to $\sigma_i^{-1}$.

Let $\Gamma_n$ denote a degree $n$ permutation chart. If there is a consistent orientation on the edges of $\Gamma_n$, then the opposite orientation will also be consistent; we let  $\Gamma^\pm_n$  denote the resulting braid charts. As before, let $M^2(\Gamma)$ denote the $n$-fold irregular simple branched cover of $D^2$ that is associated to the chart $\Gamma = \Gamma_n$, and let $\hat{M}(\Gamma)$ denote its extension to the $2$-sphere. 

\begin{theorem}
\label{orient2}
 If $\Gamma^\pm_n$ is a braid chart, then there is an embedding, $\tilde{F}: \hat{M}(\Gamma) \rightarrow S^2 \times [0,1]\times[0,1]$ of the $n$-fold irregular simple branched cover $F:\hat{M}(\Gamma) \rightarrow S^2$ such that the composition $p \circ \tilde{F}$ agrees with the covering map $F$ where $p: S^2 \times [0,1]\times [0,1] \rightarrow S^2$ is the projection onto the first factor. In this way the covering $F$ has a folded embedding.  \end{theorem}
{\sc Proof.} The lifting of each permutation to a braid induces an embedding of the surface $\hat{M}$ into $S^2 \times [0,1]\times[0,1]$. The projection to the first two factors $S^2 \times [0,1]$ induces the  generic map of Proposition~\ref{perm}. $\Box$

\subsection{A partial lifting}

It is possible that a given chart cannot be oriented consistently so that all white vertices have a flow. 
A {\it semi-oriented chart} is a permutation chart $\Gamma^*$ that includes a fourth type of vertex which is bivalent, and the chart is oriented such that 
(1) each bivalent vertex is either a  source or a  sink 
and (2) each white vertex  and each crossing has a flow.  
At a {\it source} vertex the two  emanating  edges have the same label   and point away from the vertex. At a {\it sink} the edges point towards the vertex. 

\begin{lemma} Any permutation chart can be semi-oriented.   
\end{lemma}
{\sc Proof.} Locally orient the edges so that there is a flow through all white vertices and all crossings. If there is an edge whose endpoints are either white vertices or crossings
such that the local orientations do not match, then introduce a source or sink. Do the same for all such edges.
$\Box$

\begin{lemma} A semi-oriented chart $\Gamma^*=\Gamma^*_n$ induces an immersion $\tilde{F}: \hat{M}^2(\Gamma^*) \looparrowright S^2 \times [0,1] \times [0,1]$ of the associated irregular simple branched covering $F: \hat{M}^2(\Gamma^*) \rightarrow S^2$ such that $p \circ \tilde{F} =F$ where $p$ is the projection onto the first factor. In this way the covering $F$ has a folded immersion. \end{lemma} 

Note that this result gives Theorem~\ref{mainimmersed} for $k=0$. Moreover, we may assume that the chart is standard. Clearly, a standard chart can be oriented. So by Theorem~\ref{orient2}, we have a stronger result for surfaces.

\begin{theorem} The $n$-fold simple irregular branched cover of $S^2$ branched over a finite set (even number) of points has an embedded folding. \end{theorem}

\section{Two-fold branched coverings of $S^3$ branched along a knot}

\begin{wrapfigure}[17]{l}{3.6in}
 \vspace{-.3in}
 \vspace{0.2in}
\includegraphics[width=3.5in]{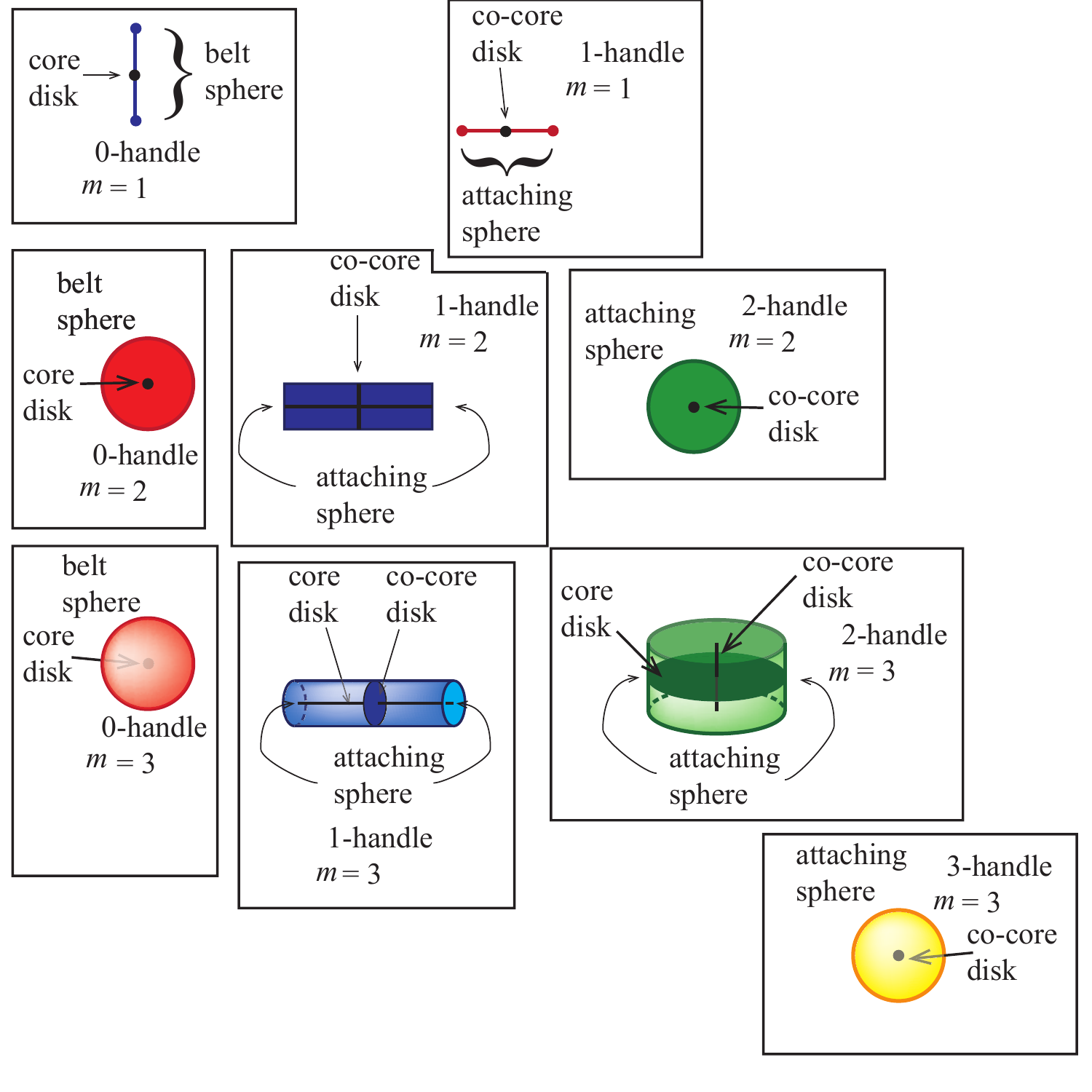}
\end{wrapfigure} Our construction of the  folding of the two-fold branched cover of $S^3$ branched along a knot or link will give rise to a Heegaard diagram for this $3$-manifold in an interesting and natural fashion. In subsequent sections, we also see handle decompositions for both $3$ and $4$ dimensional branched coverings. So for the readers convenience, we review the basic handle terminology.

\subsection{Review of handles}

Let $m$ be a fixed positive integer. Here we are mostly interested in $m=0,1,2,3$ or $4$. Let $j \in \{0,1, \ldots, m\}$.
A {\it j-handle} in an $m$-manifold $N$ is a subset homeomorphic to an $m$-disk, but decomposed as the Cartesian product $D^j \times D^{m-j}$. The subset $S^{j-1}\times D^{m-j}$ of the boundary is called the {\it attaching region} or {\it $A$-region}. The sphere $S^{j-1} \times \{0\}$ is called the {\it attaching sphere ($A$-sphere)}. The $A$-sphere is the boundary of the {\it core disk} $D^{j}\times  \{0\}$. The {\it belt region ($B$-region)} is the subset $D^j \times S^{m-j-1}$ of the boundary. 
 The {\it belt sphere ($B$-sphere)} is $\{ 0 \} \times S^{m-j-1}$. It is the boundary of the {\it co-core disk} $\{0\} \times D^{m-j}$.

\subsection{The Seifert surface as a $2$-dimensional chart}

\begin{wrapfigure}{l}{3in}\vspace{-.2in}
\includegraphics[scale=.5]{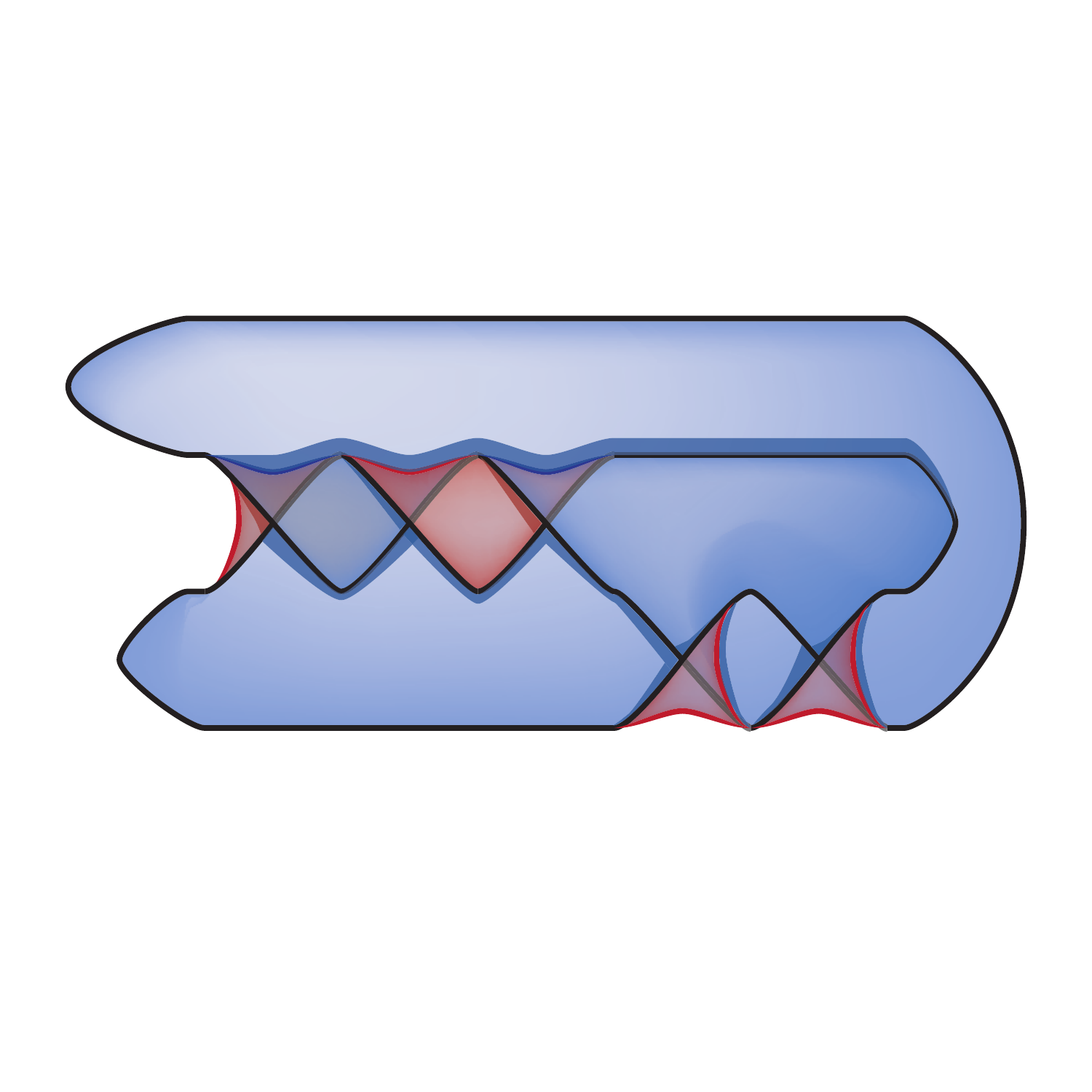}\vspace{-.2in}
\end{wrapfigure} 
There is a well-known construction for the $2$-fold branched cover of $S^3$ branched along a knot or link. First, we choose an orientable Seifert surface for the knot or link. Then we take $2$-copies of $S^3$, cut each along the Seifert surface, and glue the positive side of one copy to the negative side of the other. Our construction uses the Seifert surface as a $2$-dimensional chart. Specifically, we choose a height function on $S^3$ that restricts to a non-degenerate height function on the Seifert surface of the knot or link.

By convention we arrange the knot and the Seifert surface with the initial minimum to the left and the final maximum to the right. The knot may be taken to lie within a $3$-disk that is parametrized as $[0,1]\times [0,1]\times [0,1]$. The first factor is thought to run left-to-right, the second from bottom-to-top, and the third fact is back-to-front. 
We often rotate each still 90 degrees in clockwise (or anti-clockwise) direction.   Then the second factor is thought to run left-to-right. 
The projection onto the first factor is a non-degenerate height function for the Seifert surface. 
\begin{wrapfigure}[17]{l}{4.25in}\vspace{-.4in}
\includegraphics[scale=.7]{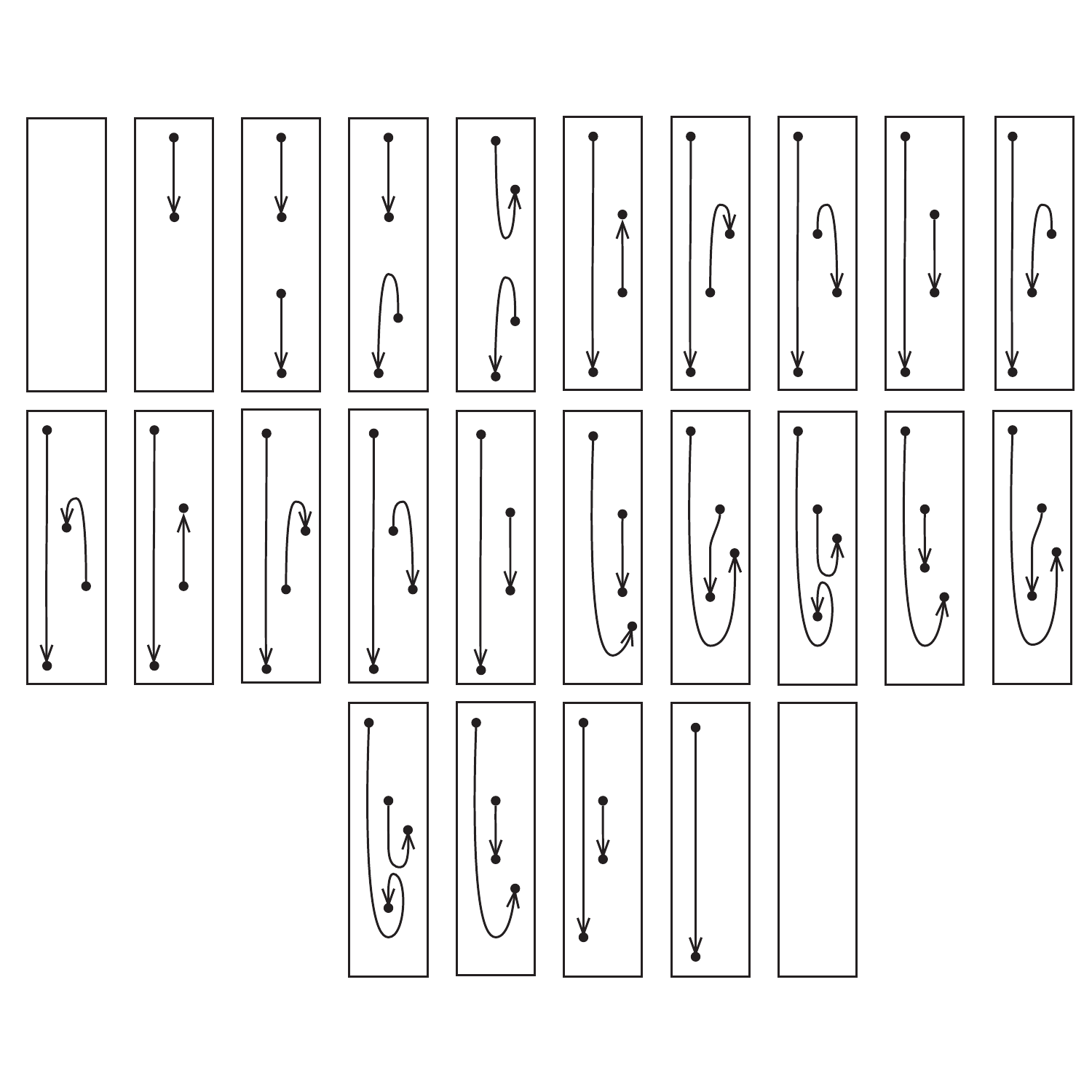}\vspace{-.4in}
\end{wrapfigure}
In particular each critical point for the surface is at a different level, each crossing is at a different level,  the surface has optimal points at the optimal points of the link (it may have other optima), and it may have saddle points. We cut the link and the Seifert surface between critical levels where the crossings are considered to be critical. Since Seifert surfaces are constructed by twisted bands, we will also arrange that the twists are achieved by particular critical motions near the crossings. The illustration above and to the left indicate these ideas for the knot $5_2$.

For a Seifert surface $\mathcal{F} \subset [0,1]\times [0,1]\times [0,1]$ and for each $s \in [0,1]$, let 
$\Gamma_s = \Gamma_s(\mathcal{F}) \subset [0,1]\times [0,1]$ denote the cross-section of $\mathcal{F}$ at $s$, i.e., 
$$
\mathcal{F} \cap \{s\} \times [0,1]\times [0,1] = \{s\} \times \Gamma_s. 
$$
For a non-critical level $s$, the cross-section $\Gamma_s$ is 
a permutation chart of degree $2$. An orientation of the link induces an orientation of the Seifert surface, and an orientation for the cross-sectional charts. The conventions are indicated 
in the figure 
above, 
when we rotate the stills 90 degrees: an arc pointing towards the right is positive and the corresponding black vertex is a sink; a left-pointing arc is negative and the corresponding black vertex is a source. 

For each non-critical cross-section, 
the oriented cross-sectional chart $\Gamma_s$ is a braid chart of degree $2$ and it describes an embedding 
$F_s \hookrightarrow S^2 \times D^2$ of a closed oriented surface $F_s$ as a surface braid (cf. \cite{Kam2002}).  
We will call $F_s \subset S^2 \times D^2$ the {\it cross-sectional surface} or the {\it cross-sectional surface braid} at $s$.  

The rest of the construction is to describe the changes to surface braids induced by critical events for the Seifert surface. Thus we will have an embedded folded solid over the cube $[0,1]\times [0,1] \times [0,1]$, and we will extend this to the trivial cover of the complementary disk in $S^3$.

\subsection{Illustrations of the critical events}

Here is a list of the critical events.

\noindent 
 {\bf 1.} $\stackrel{1-H}{\Longrightarrow}$ or $\stackrel{2-H}{\Longrightarrow}$  ---{\it handle attachment}. A local optimal  point of knot diagram occurs. The move   \begin{wrapfigure}[10]{r}{3in}\vspace{-.3in}
\includegraphics[scale=.4]{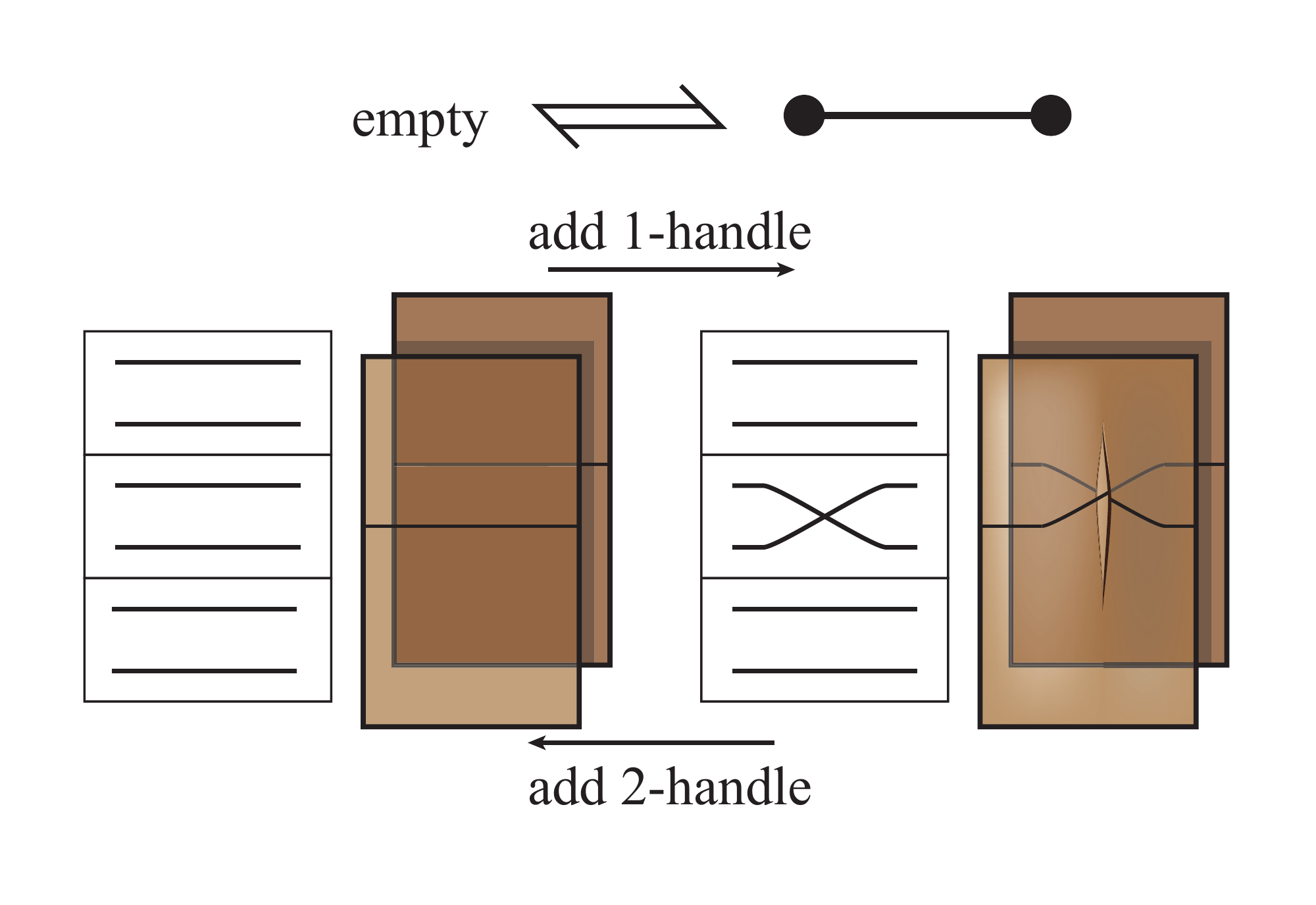} 
\end{wrapfigure} from left-to-right 
in the illustration to the right
represents a $1$-handle being attached between the cross-sectional surfaces. The moves from right-to-left represents a $2$-handle being attached between the cross-sectional surfaces. The attaching region for the $1$-handle consist of a pair of disks on the left-hand-side of the figure. The belt region is an annular neighborhood of the double curve. It is twisted since the normal orientations of the two disks involved are parallel (as opposed to anti-parallel). 

\noindent 
\begin{wrapfigure}{l}{3in}\vspace{-0.2in}
\includegraphics[width=2.5in]{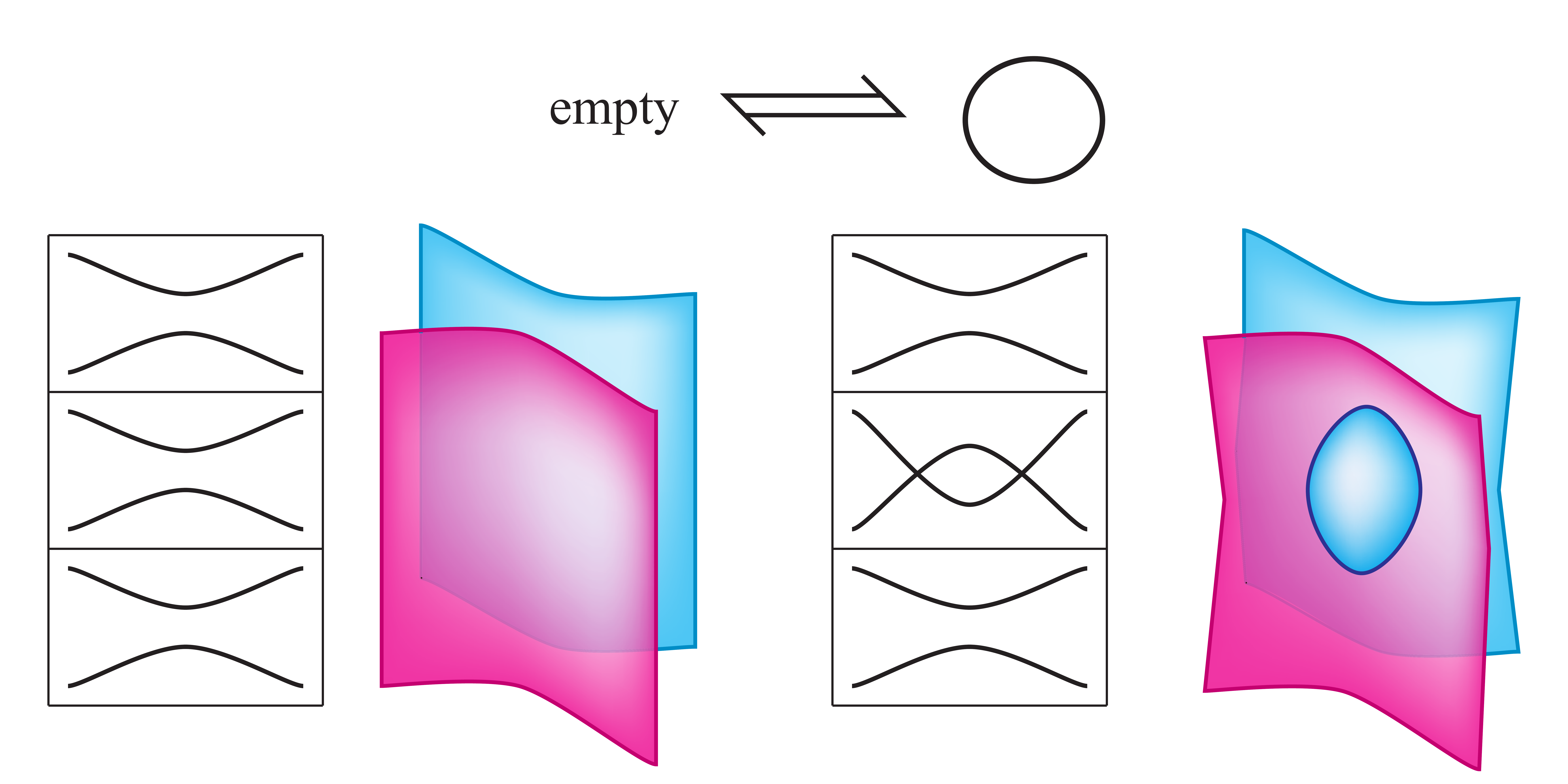} 
\end{wrapfigure}
{\bf 2.}
 $\stackrel{IIb}{\Longrightarrow}$  --- {\it a type II bubble move}.  An optimum occurs in the interior of the Seifert surface. In this case the successive cross-sections differ by the inclusion of a simple closed curve. The cross-sectional surfaces are isotopic as surface braids since they differ by  the chart move that is indicated.

\begin{wrapfigure}{r}{3in}\vspace{-0.2in}
\includegraphics[width=2.5in]{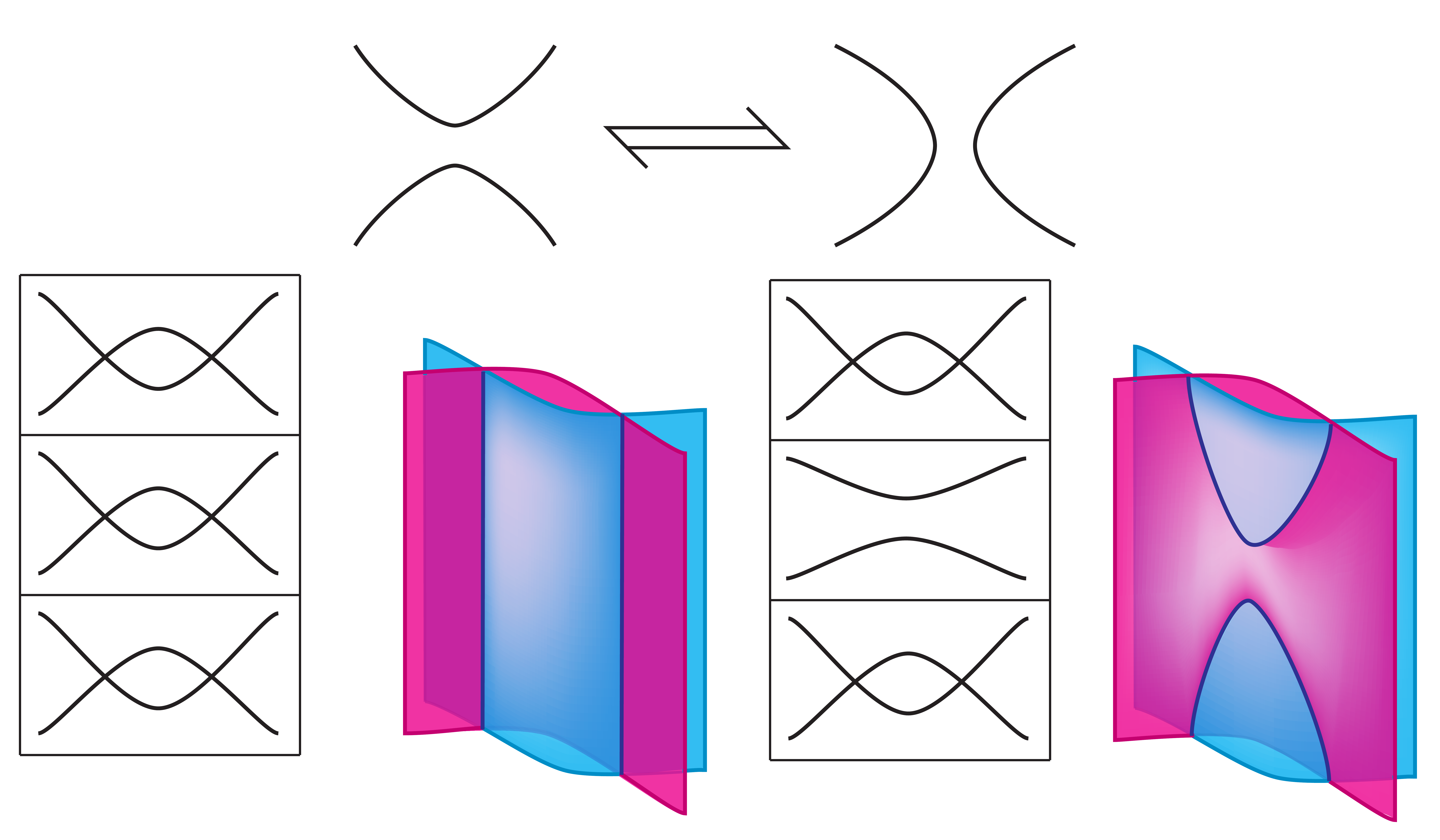} 
\end{wrapfigure}
\noindent {\bf 3.} $\stackrel{IIs}{\Longrightarrow}$ --- {\it A type II saddle move}. A saddle occurs in the interior of the Seifert surface. In this case, the relation $\tau_1 \cdot \tau_1 \rightarrow 1$ is followed by the relation $1 \rightarrow \tau_1 \cdot \tau_1$ in the permutation representation, and the sequence 
$\tau_1 \cdot \tau_1 \rightarrow 1  \rightarrow \tau_1 \cdot \tau_1$ is replaced by the identity sequence $\tau_1 \cdot \tau_1 \rightarrow   \tau_1 \cdot \tau_1 \rightarrow \tau_1 \cdot \tau_1$  (or vice versa). In the braid case, the sequence $\sigma_1^{\pm 1} \cdot \sigma_1^{\mp 1} \rightarrow 1  \rightarrow \sigma_1^{\pm 1} \cdot \sigma_1^{\mp 1}$ is replaced by the identity sequence
$\sigma_1^{\pm 1} \cdot \sigma_1^{\mp 1} \rightarrow \sigma_1^{\pm 1} \cdot \sigma_1^{\mp 1}.$ 
The cross-sectional surfaces are isotopic as surface braids since they differ by  the chart move that is indicated.

 \begin{wrapfigure}[6]{l}{3 in}\vspace{-.1in}
\includegraphics[width=2.75in]{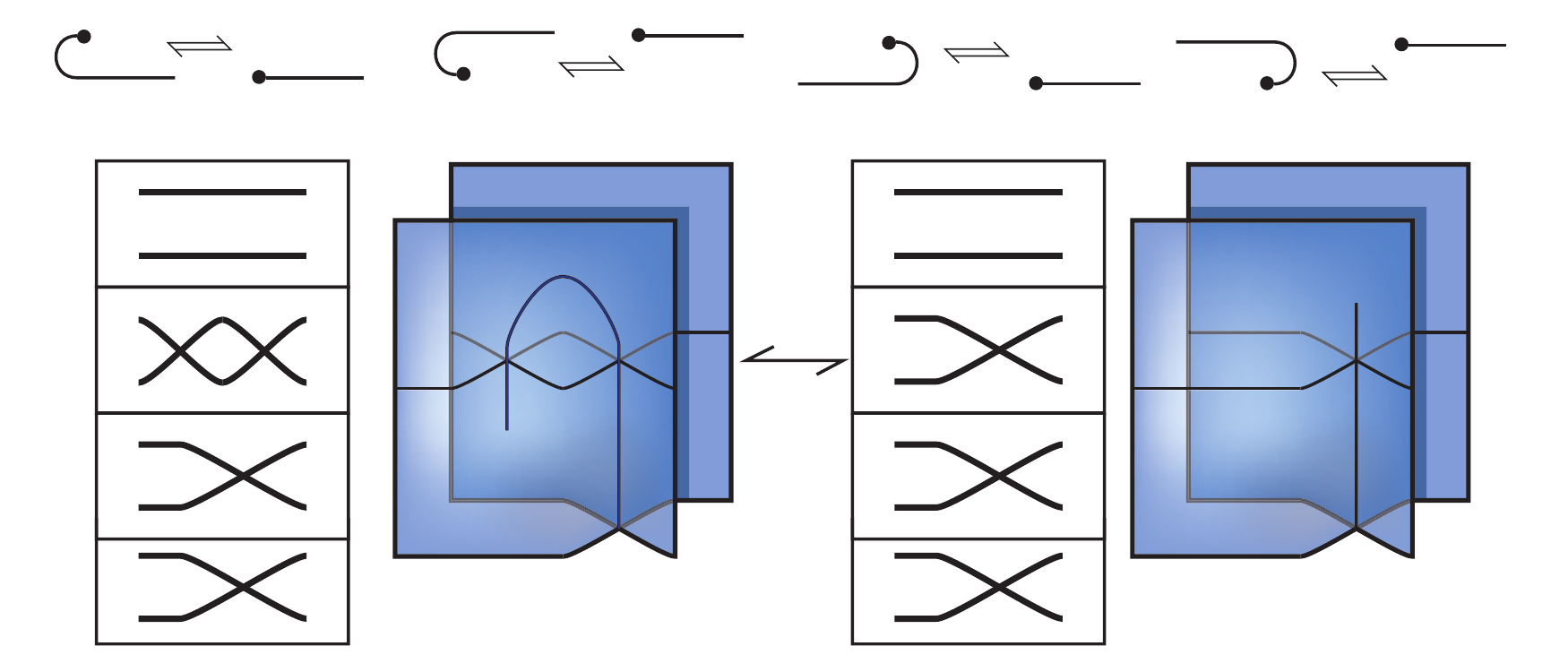} 
\end{wrapfigure}
\noindent 
{\bf 4.} 
 $\stackrel{CC}{\Longrightarrow}$ --- {\it candy-cane moves}. The Seifert surface can bend near the boundary. The possible changes in the cross-sections are depicted on the top line of the figure to the left. On the bottom, we see that a type-II move followed by a branch point can be replaced by a branch point in the opposite direction. The illustration indicates the phenomenon in the permutation case, but the braid case is also easy to understand. In this case the crossing information is consistent. 
 The cross-sectional surfaces are isotopic as surface braids, although the branch points corresponding the black vertices  move slightly.

 \begin{wrapfigure}[7]{r}{2 in}\vspace{-.1in}\vspace{-.3in}
\includegraphics[width=1.75in]{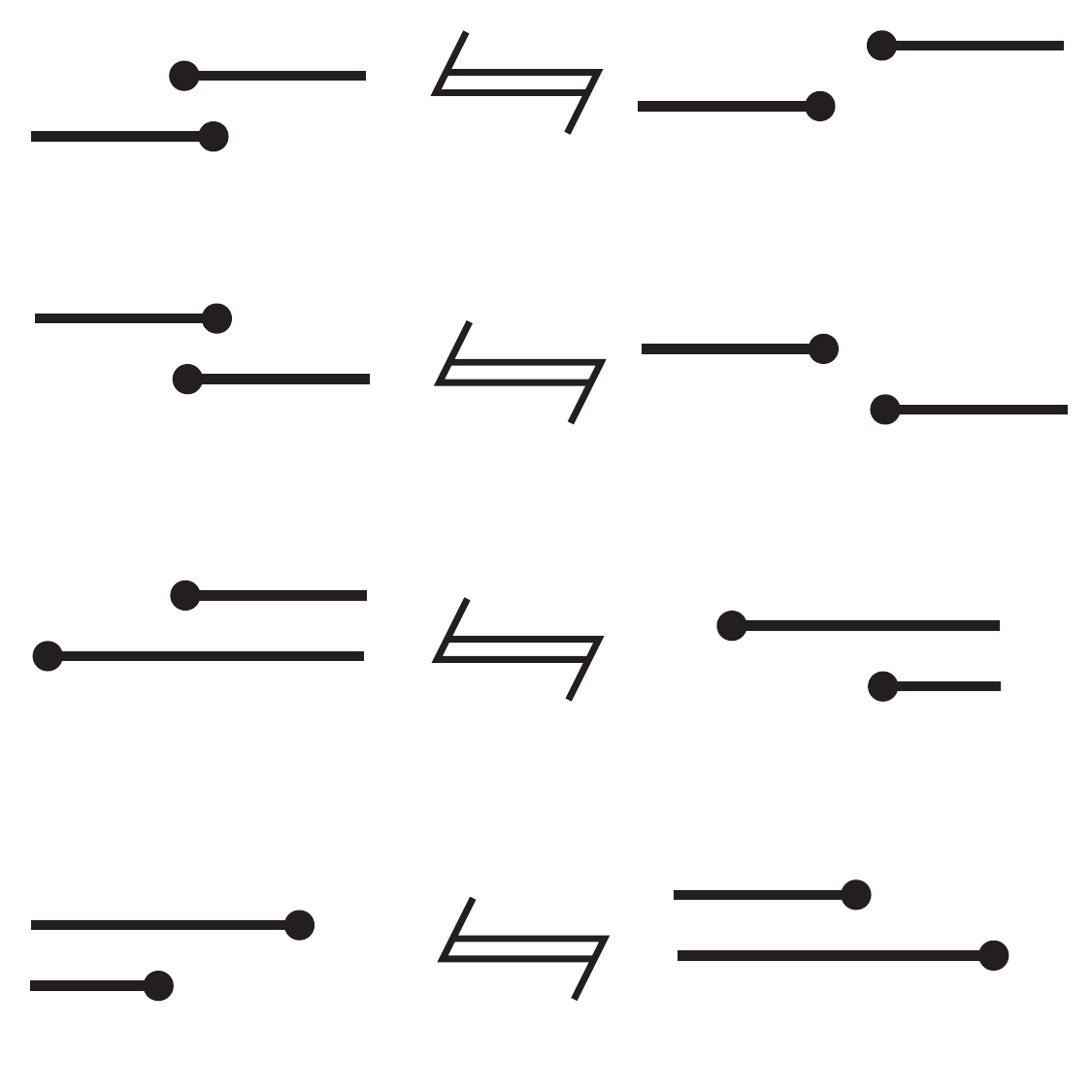} 
\end{wrapfigure}
\noindent {\bf 5.} $\stackrel{X}{\Longrightarrow}$  --- {\it exchange moves.}  
 End points of arcs exchange places. These moves occur near a crossing in the diagram. The illustration given indicates the variety of ways in which these exchanges can take place. Within the surfaces constructed from the cross-sections, the situation is quite straight-forward. The branch points and the double point arcs that terminate at the branch points are free to move within the surfaces as long as the arcs of double points do not intersect.
  The cross-sectional surfaces are isotopic as surface braids because the branch points corresponding the black vertices  move slightly.

 \begin{wrapfigure}[7]{l}{3 in}\vspace{-.1in}
\includegraphics[width=2.75in]{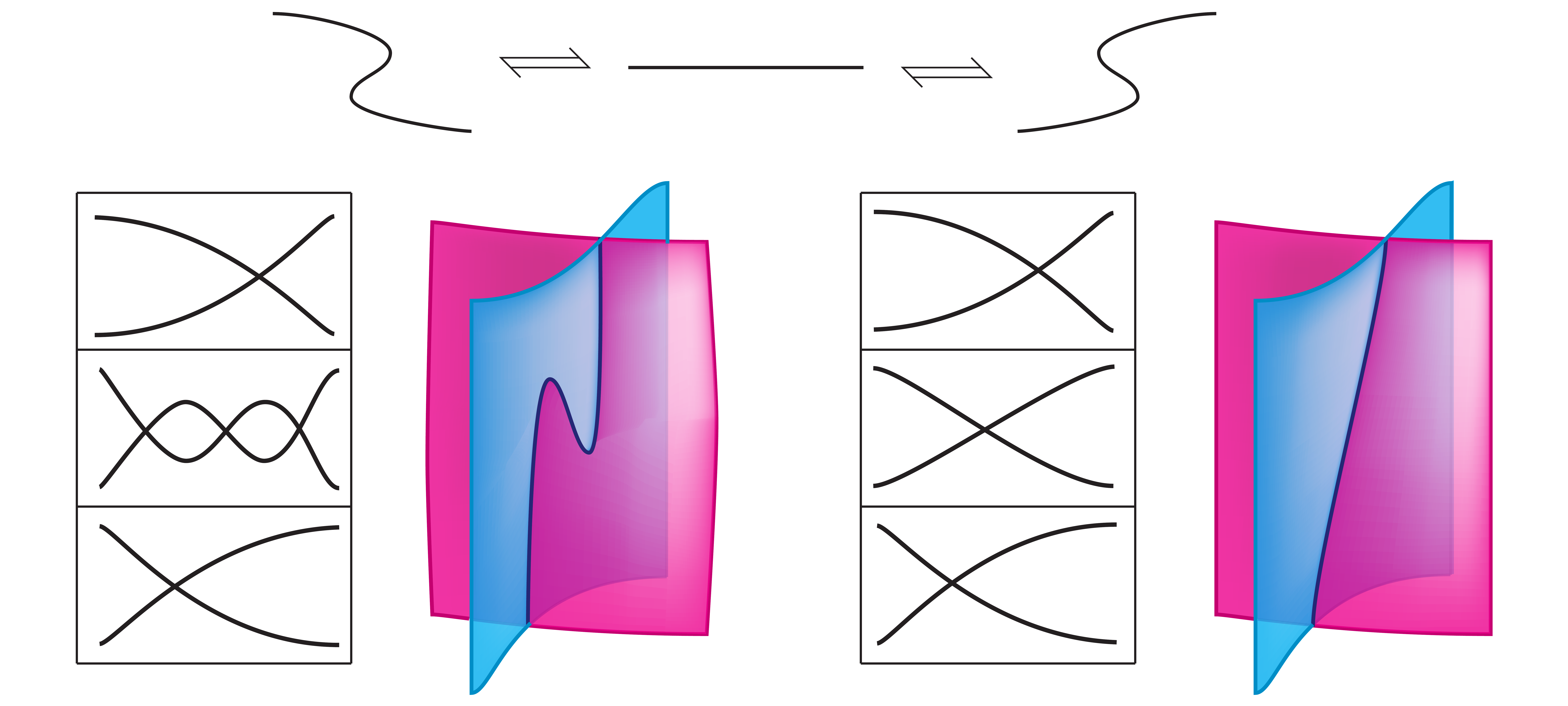} 
\end{wrapfigure}\noindent
{\bf 6.} $\stackrel{Z}{\Longrightarrow}$ --- {\it type II zig-zag moves}. Cusps occur on the Seifert surface. One such cusp is indicated below.  The cross-sectional surfaces are isotopic as surface braids because under-lying charts are topologically equivalent.

\begin{lemma} A given segment in a braid chart can be split into two via the moves above. Thus splitting an edge into two edges is also caused by adding a disjoint segment. In either case, a handle is attached to the underlying $3$-manifold:
$$\bullet \mkern-8mu \raisebox{0.04in}{ \rule{30pt}{.7pt} }\mkern-8mu \bullet  \ \Leftrightarrow  \   \bullet \mkern-10mu \raisebox{0.04in}{ \rule{32pt}{.7pt} }\mkern-12mu \bullet  \ \ \bullet \mkern-8mu \raisebox{0.04in}{ \rule{32pt}{.7pt} }\mkern-12mu \bullet . $$ 
\end{lemma}
{\sc Proof.} 
\vspace{-.5cm}
\begin{center}
\includegraphics[scale=.55]{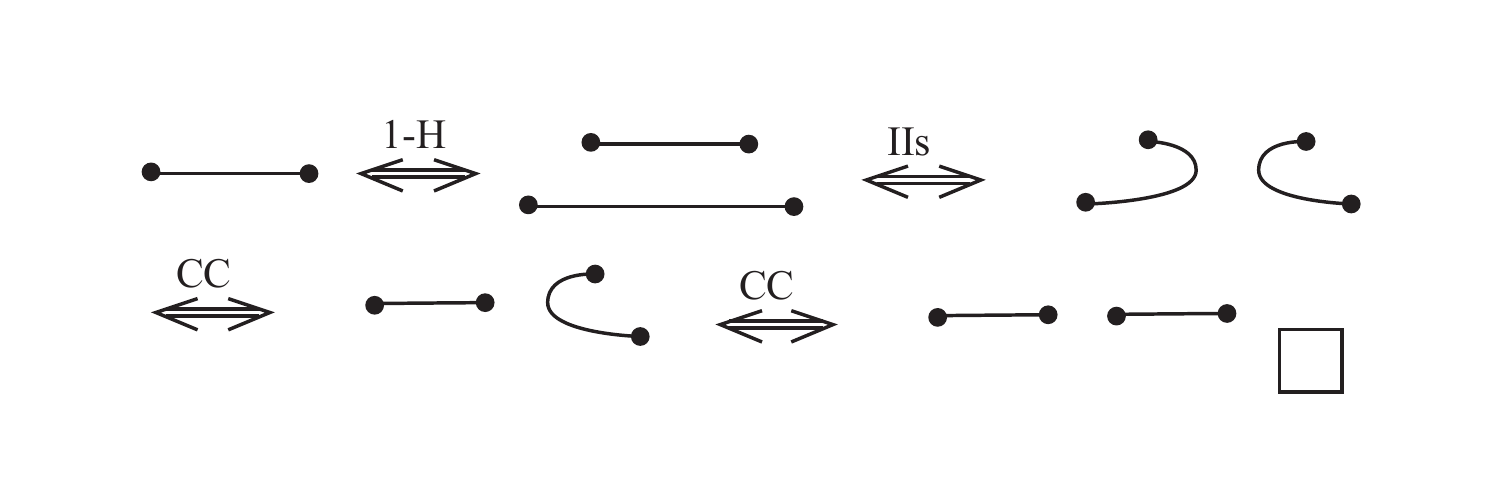}
\end{center}

\subsection{Example}

\begin{wrapfigure}{l}{3in}\vspace{-.2in}
\includegraphics[scale=.5]{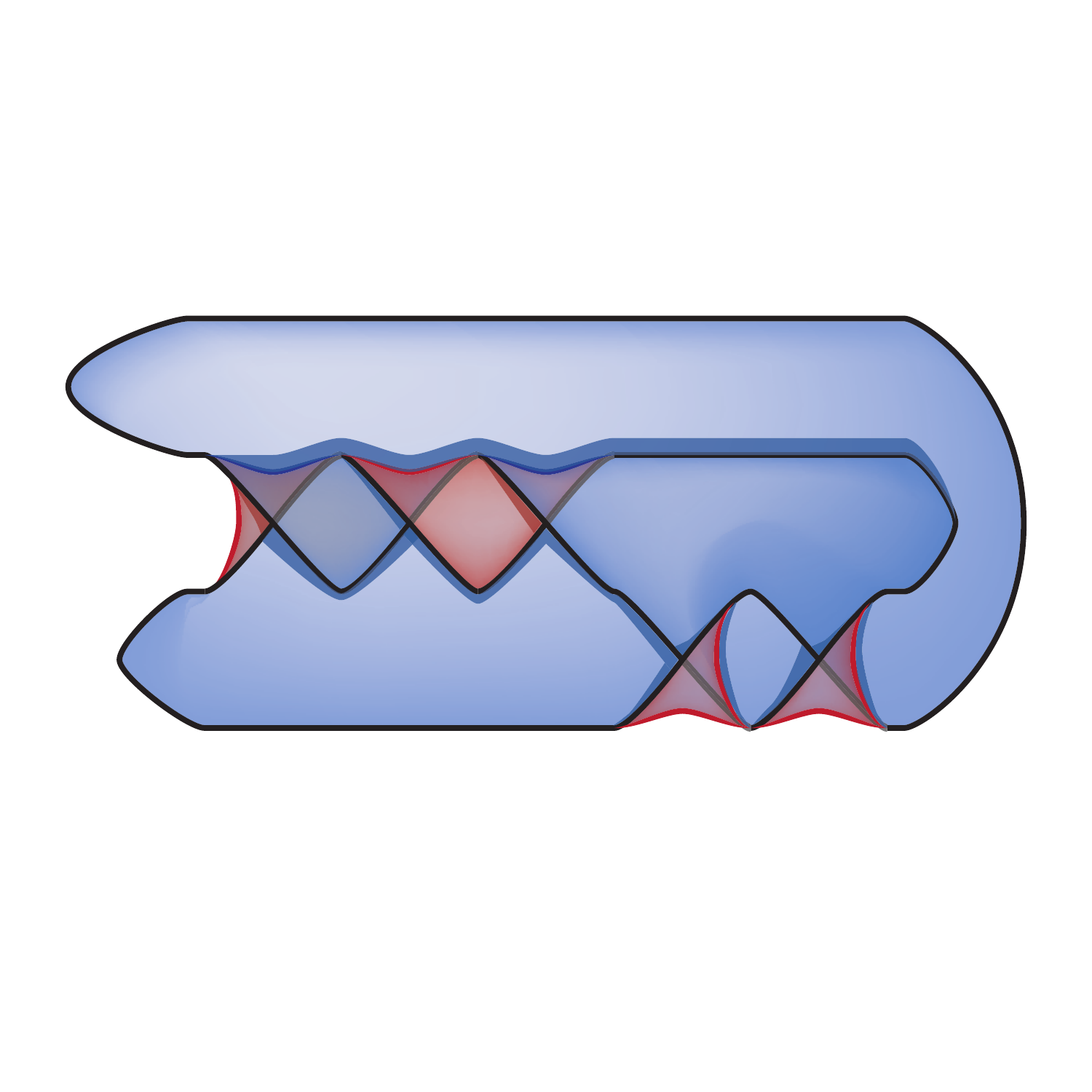}\vspace{-.2in}
\end{wrapfigure} In this section, we give an explicit description of the $2$-fold branched cover of $S^3$ branched along the knot $5_2$.

For the reader's convenience, we reproduce the illustration of the Seifert surface for $5_2$ to the left.
 Recall, that the knot is embedded in $[0,1]\times [0,1]\times [0,1]$ and projection onto the first factor provides a height function for the knot and for the Seifert surface. Thus critical points and crossings occur at distinct levels. Each generic cross section $p_1^{-1}(t_i)$ consists of a braid chart where the orientation in the chart is induced from the orientation of the knot: a right pointing arc in the knot corresponds to the target of an arrow in the chart. The charts in the cross-sections are oriented in the figure so that the right edges correspond to the portion of the knot diagram that is closest to the observer.

We read each cross-sectional chart as a sequence of braid words. The sequence always starts and ends at the empty word. A  horizontal line intersecting the chart at a non-critical level intersects the arcs of a chart. If the arc is up pointing, the braid generator is positive; if down-pointing, then a negative generator is encountered. The critical points of the chart $\cap$ and $\cup$ correspond to type II moves. The black vertices correspond to branch points --- the introduction or deletion of a braid generator or its inverse. 

A critical event between charts is one of the changes in charts catalogued in the previous section. 

Thus we have \begin{itemize}
\item for each cross section of a chart a braid word;
\item for the set of cross-sections a ``paragraph of braid words;"
\item for the collection of charts an ``essay" that consists of paragraphs of braid words;
\item methods of getting from one word to another;
\item methods of getting from one paragraph to the next;
\item some standard introductory words and phrases that lead to the paragraphs and correspond to empty charts.
\end{itemize}

\begin{wrapfigure}{r}{2.5in}\vspace{-.35in}
\includegraphics[scale=.15]{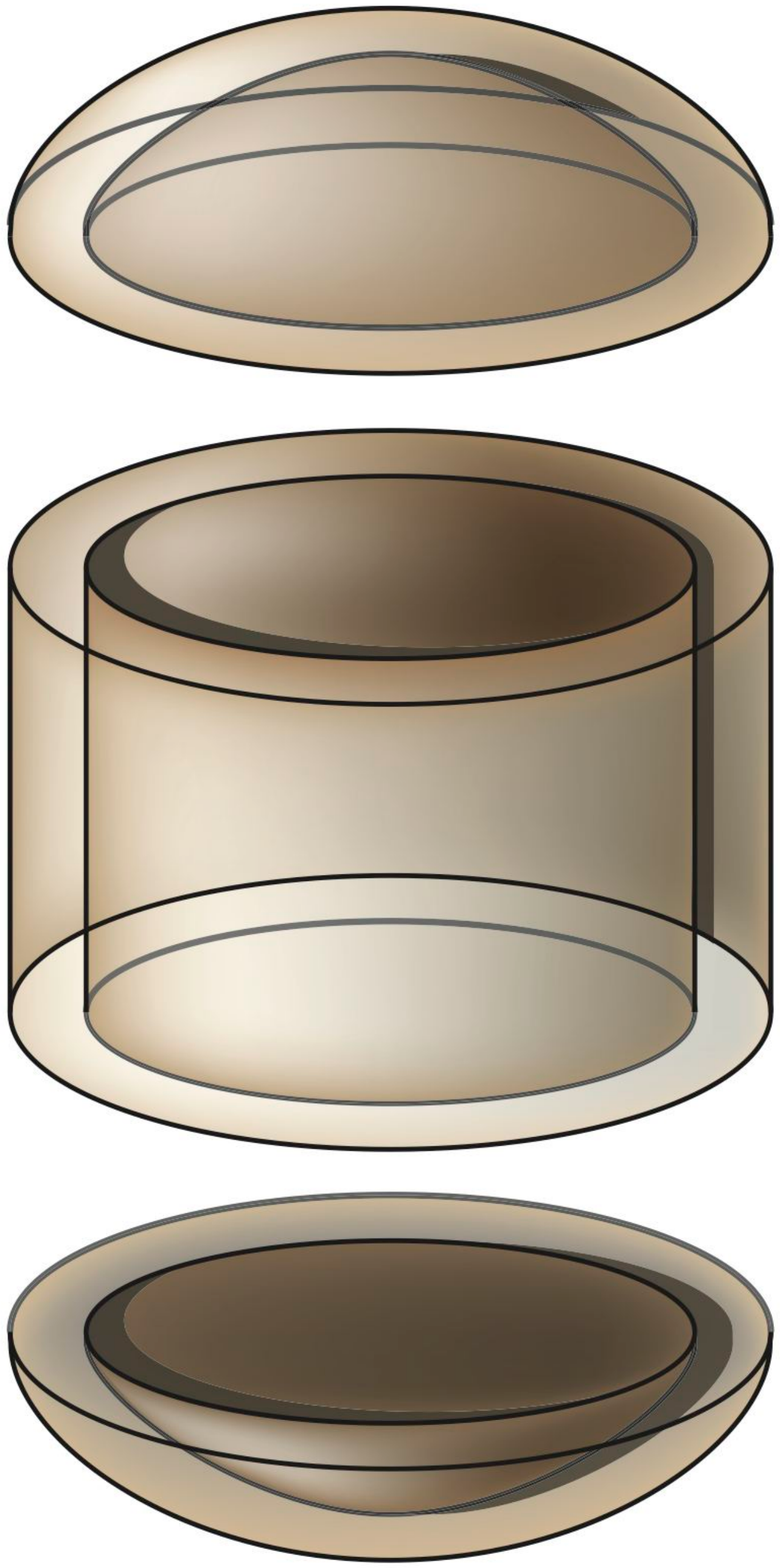}\end{wrapfigure}
Before developing the braid word essay associated to $5_2$, we  discuss the introductory words and paragraphs.
Each braid word represents a two-fold  cover of a circle. The representation is achieved by taking the braid closure of each word. In particular, the braid closure of the empty word corresponds to two nested circles in the plane. Since the chart is embedded as a $2$-disk in the $2$-sphere,  the unknotted, nested circles at the top and the bottom sequentially bound disks that cover the polar caps of the $2$-spheres. The  closure of an empty chart is depicted to the right. More generally, the pair of nested caps at the top of the diagram  and the pair of nest bowls at the bottom appear at the top and bottom of any particular chart.

Thus the empty charts that appear in the cross-sections to the immediate left and right of the knot with Seifert surface represent a pair of nested, embedded spheres. These spheres, in turn, bound a pair embedded $3$-balls  (at each end of the chart sequence) that trivially cover the left and right ends of $S^3$. By turning the picture of the embedded knot on its side, we can imagine the $3$-sphere as a northern polar cap which is a $3$-ball,  a temperate zone which is homeomorphic to $S^2\times [0,1]$ and in which the knot lies, and a southern polar cap which is another $3$-ball. 

Returning to the specifics of the sequence of charts for the Seifert surface of the knot $5_2$, we have the following. We start with a pair of nested $2$-spheres that were created as in the preceding paragraph and represented as an empty chart. At the end of each of the sequences below, we have a pair of nested $2$-spheres that successively bound a pair of $3$-balls. The notation ${\mbox{\rm b}}_j^{\pm}$ indicates that a crossing is added or subtracted via a black vertex at the $j\/$th position in the word. Thus a segment is flanked by ${\mbox{\rm b}}_j^{+}$  and ${\mbox{\rm b}}_{j+1}^{-}$.
The notation ${\mbox{\rm II}}_j^{\pm}$ indicates that a type II move has been performed that either adds or subtracts a pair of oppositely signed braid generators with the first insertion or deletion occurring at the $j\/$th position of the next word.

\noindent
{\bf The braid essay associated to the Seifert surface of the knot $5_2$:}

\vspace{.5cm}

\noindent
\begin{wrapfigure}[7]{l}{2.9in}\vspace{-.4in}\vspace{-.2in}
\includegraphics[scale=.42]{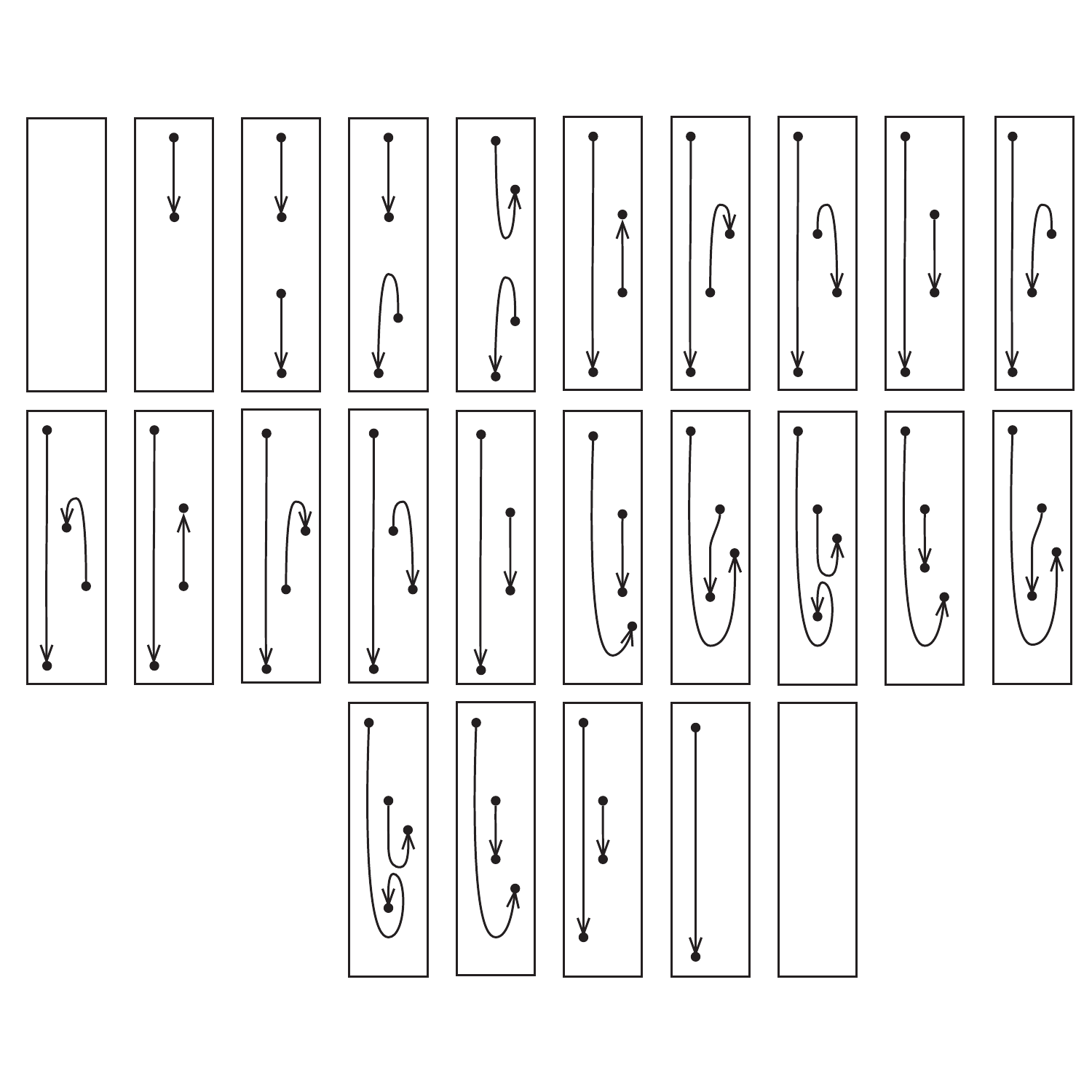}\vspace{-.4in}
\end{wrapfigure}
 $[(\emptyset)]$ $\stackrel{1-H}{\Longrightarrow}$
 
 \noindent
 $[(\emptyset) \stackrel{\mbox{\rm b}_0^+}{\rightarrow} 
 (\sigma_1^{-1}) \stackrel{\mbox{\rm b}_1^-}{\rightarrow}
(\emptyset)]$
 $\stackrel{1-H}{\Longrightarrow}$
 
  \noindent
 $[(\emptyset) \stackrel{\mbox{\rm b}_0^+}{\rightarrow} 
 (\sigma_1^{-1}) \stackrel{\mbox{\rm b}_1^-}{\rightarrow}
(\emptyset) \stackrel{\mbox{\rm b}_0^+}{\rightarrow} 
 (\sigma_1^{-1}) \stackrel{\mbox{\rm b}_1^-}{\rightarrow}
(\emptyset) 
]$
 $\stackrel{CC}{\Longrightarrow}$

  \noindent
 $[(\emptyset) \stackrel{\mbox{\rm b}_0^+}{\rightarrow} 
 (\sigma_1^{-1}) \stackrel{\mbox{\rm b}_1^-}{\rightarrow}
(\emptyset) \stackrel{\mbox{\rm II}_1^+}{\rightarrow} 
 (\sigma_1^{-1}\sigma_1)
\stackrel{\mbox{\rm b}^-_2}{\rightarrow} 
 (\sigma_1^{-1}) \stackrel{\mbox{\rm b}_1^-}{\rightarrow}
(\emptyset) 
]$
 $\stackrel{CC}{\Longrightarrow}$

  \noindent
 $[(\emptyset) \stackrel{\mbox{\rm b}_0^+}{\rightarrow} 
 (\sigma_1^{-1}) \stackrel{\mbox{\rm b}_1^+}{\rightarrow}
 (\sigma_1^{-1}\sigma_1) \stackrel{\mbox{\rm II}_1^-}{\rightarrow}
(\emptyset) \stackrel{\mbox{\rm II}_1^+}{\rightarrow} 
 (\sigma_1^{-1}\sigma_1)
\stackrel{\mbox{\rm b}^-_2}{\rightarrow} 
 (\sigma_1^{-1}) \stackrel{\mbox{\rm b}^-_1}{\rightarrow}
(\emptyset) 
]$
 $\stackrel{IIs}{\Longrightarrow}$
 
  \noindent
 $[(\emptyset) \stackrel{\mbox{\rm b}_0^+}{\rightarrow} 
 (\sigma_1^{-1}) \stackrel{\mbox{\rm b}_1^+}{\rightarrow}
 (\sigma_1^{-1}\sigma_1) 
\stackrel{\mbox{\rm b}^-_2}{\rightarrow} 
 (\sigma_1^{-1}) \stackrel{\mbox{\rm b}^-_1}{\rightarrow}
(\emptyset) 
]$
 $\stackrel{CC}{\Longrightarrow}$

   \noindent
 $[(\emptyset) \stackrel{\mbox{\rm b}_0^+}{\rightarrow} 
 (\sigma_1^{-1}) \stackrel{\mbox{\rm II}_2^+}{\rightarrow}
 (\sigma_1^{-1}\sigma_1\sigma_1^{-1}) 
\stackrel{\mbox{\rm b}^-_3}{\rightarrow} 
 (\sigma_1^{-1}\sigma_1) \stackrel{\mbox{\rm b}^-_2}{\rightarrow} (\sigma_1^{-1}) \stackrel{\mbox{\rm b}_1^-}{\rightarrow}
(\emptyset) 
]$
 $\stackrel{X}{\Longrightarrow}$

  \noindent
 $[(\emptyset) \stackrel{\mbox{\rm b}_0^+}{\rightarrow} 
 (\sigma_1^{-1}) \stackrel{\mbox{\rm II}_2^+}{\rightarrow}
 (\sigma_1^{-1}\sigma_1\sigma_1^{-1}) 
\stackrel{\mbox{\rm b}^-_2}{\rightarrow} 
 (\sigma_1^{-1}\sigma_1^{-1}) \stackrel{\mbox{\rm b}^-_2}{\rightarrow} (\sigma_1^{-1}) \stackrel{\mbox{\rm b}^-_1}{\rightarrow}
(\emptyset) 
]$
 $\stackrel{CC}{\Longrightarrow}$
 
   \noindent
 $[(\emptyset) \stackrel{\mbox{\rm b}_0^+}{\rightarrow} 
 (\sigma_1^{-1}) \stackrel{\mbox{\rm b}_1^+}{\rightarrow}
 (\sigma_1^{-1}\sigma_1^{-1}) 
\stackrel{\mbox{\rm b}^-_2}{\rightarrow} 
 (\sigma_1^{-1}) \stackrel{\mbox{\rm b}^-_1}{\rightarrow}
(\emptyset) 
]$
 $\stackrel{CC}{\Longrightarrow}$

   \noindent
 $[(\emptyset) \stackrel{\mbox{\rm b}_0^+}{\rightarrow} 
 (\sigma_1^{-1}) \stackrel{\mbox{\rm II}_2^+}{\rightarrow}
 (\sigma_1^{-1}\sigma_1^{-1}\sigma_1) 
\stackrel{\mbox{\rm b}^-_3}{\rightarrow} 
 (\sigma_1^{-1}\sigma_1^{-1}) \stackrel{\mbox{\rm b}^-_2}{\rightarrow} (\sigma_1^{-1}) \stackrel{\mbox{\rm b}^-_1}{\rightarrow}
(\emptyset) 
]$
 $\stackrel{X}{\Longrightarrow}$

  \noindent
 $[(\emptyset) \stackrel{\mbox{\rm b}_0^+}{\rightarrow} 
 (\sigma_1^{-1}) \stackrel{\mbox{\rm II}_2^+}{\rightarrow}
 (\sigma_1^{-1}\sigma_1^{-1}\sigma_1) 
\stackrel{\mbox{\rm b}^-_2}{\rightarrow} 
 (\sigma_1^{-1}\sigma_1) \stackrel{\mbox{\rm b}^-_2}{\rightarrow} (\sigma_1^{-1}) \stackrel{\mbox{\rm b}_1^-}{\rightarrow}
(\emptyset) 
]$
 $\stackrel{CC}{\Longrightarrow}$
 
 \noindent
 $[(\emptyset) \stackrel{\mbox{\rm b}_0^+}{\rightarrow} 
 (\sigma_1^{-1}) 
\stackrel{\mbox{\rm b}^+_1}{\rightarrow} 
 (\sigma_1^{-1}\sigma_1) \stackrel{\mbox{\rm b}^-_2}{\rightarrow} (\sigma_1^{-1}) \stackrel{\mbox{\rm b}^-_1}{\rightarrow}
(\emptyset) 
]$
 $\stackrel{CC}{\Longrightarrow}$

   \noindent
 $[(\emptyset) \stackrel{\mbox{\rm b}_0^+}{\rightarrow} 
 (\sigma_1^{-1}) \stackrel{\mbox{\rm II}_2^+}{\rightarrow}
 (\sigma_1^{-1}\sigma_1\sigma_1^{-1}) 
\stackrel{\mbox{\rm b}^-_3}{\rightarrow} 
 (\sigma_1^{-1}\sigma_1) \stackrel{\mbox{\rm b}^-_2}{\rightarrow} (\sigma_1^{-1}) \stackrel{\mbox{\rm b}_1^-}{\rightarrow}
(\emptyset) 
]$
 $\stackrel{X}{\Longrightarrow}$

  \noindent
 $[(\emptyset) \stackrel{\mbox{\rm b}_0^+}{\rightarrow} 
 (\sigma_1^{-1}) \stackrel{\mbox{\rm II}_2^+}{\rightarrow}
 (\sigma_1^{-1}\sigma_1\sigma_1^{-1}) 
\stackrel{\mbox{\rm b}^-_2}{\rightarrow} 
 (\sigma_1^{-1}\sigma_1^{-1}) \stackrel{\mbox{\rm b}^-_2}{\rightarrow} (\sigma_1^{-1}) \stackrel{\mbox{\rm b}^-_1}{\rightarrow}
(\emptyset) 
]$
 $\stackrel{CC}{\Longrightarrow}$
 
   \noindent
 $[(\emptyset) \stackrel{\mbox{\rm b}_0^+}{\rightarrow} 
 (\sigma_1^{-1}) \stackrel{\mbox{\rm b}_1^+}{\rightarrow}
 (\sigma_1^{-1}\sigma_1^{-1}) 
\stackrel{\mbox{\rm b}^-_2}{\rightarrow} 
 (\sigma_1^{-1}) \stackrel{\mbox{\rm b}^-_1}{\rightarrow}
(\emptyset) 
]$
 $\stackrel{CC}{\Longrightarrow}$

  \noindent
 $[(\emptyset) \stackrel{\mbox{\rm b}_0^+}{\rightarrow} 
 (\sigma_1^{-1}) \stackrel{\mbox{\rm b}_1^+}{\rightarrow}
 (\sigma_1^{-1}\sigma_1^{-1}) 
  \stackrel{\mbox{\rm b}^-_2}{\rightarrow}
 (\sigma_1^{-1}) 
 \stackrel{\mbox{\rm b}^+_1}{\rightarrow}
 (\sigma_1^{-1}\sigma_1) 
 \stackrel{\mbox{\rm II}^-_1}{\rightarrow} 
(\emptyset) 
]$
 $\stackrel{X}{\Longrightarrow}$
 
  \noindent
 $[(\emptyset) \stackrel{\mbox{\rm b}_0^+}{\rightarrow} 
 (\sigma_1^{-1}) \stackrel{\mbox{\rm b}_1^+}{\rightarrow}
 (\sigma_1^{-1}\sigma_1^{-1}) 
  \stackrel{\mbox{\rm b}^+_2}{\rightarrow}
 (\sigma_1^{-1}\sigma_1^{-1}\sigma_1) 
 \stackrel{\mbox{\rm b}_2^-}{\rightarrow}
 (\sigma_1^{-1}\sigma_1) 
 \stackrel{\mbox{\rm II}^-_1}{\rightarrow} 
(\emptyset) 
]$
 $\stackrel{IIs}{\Longrightarrow}$

  \noindent
 $[(\emptyset) \stackrel{\mbox{\rm b}_0^+}{\rightarrow} 
 (\sigma_1^{-1}) \stackrel{\mbox{\rm b}_1^+}{\rightarrow}
 (\sigma_1^{-1}\sigma_1^{-1}) 
  \stackrel{\mbox{\rm b}^+_2}{\rightarrow}
 (\sigma_1^{-1}\sigma_1^{-1}\sigma_1)
  \stackrel{\mbox{\rm II}^-_2}{\rightarrow}
  (\sigma_1^{-1})
  \stackrel{\mbox{\rm II}^+_2}{\rightarrow}
 (\sigma_1^{-1}\sigma_1^{-1}\sigma_1)
  \stackrel{\mbox{\rm b}_2^-}{\rightarrow}
 (\sigma_1^{-1}\sigma_1) 
 \stackrel{\mbox{\rm II}^-_1}{\rightarrow} 
(\emptyset) 
]$
 $\stackrel{CC^2}{\Longrightarrow}$  
 
   \noindent
 $[(\emptyset) \stackrel{\mbox{\rm b}_0^+}{\rightarrow} 
 (\sigma_1^{-1}) \stackrel{\mbox{\rm b}_1^+}{\rightarrow}
 (\sigma_1^{-1}\sigma_1^{-1}) 
  \stackrel{\mbox{\rm b}_2^-}{\rightarrow}
 (\sigma_1^{-1})
  \stackrel{\mbox{\rm b}_1^+}{\rightarrow}
 (\sigma_1^{-1}\sigma_1) 
 \stackrel{\mbox{\rm II}_1^-}{\rightarrow} 
(\emptyset) 
]$
 $\stackrel{X}{\Longrightarrow}$
 
   \noindent
 $[(\emptyset) \stackrel{\mbox{\rm b}_0^+}{\rightarrow} 
 (\sigma_1^{-1}) \stackrel{\mbox{\rm b}_1^+}{\rightarrow}
 (\sigma_1^{-1}\sigma_1^{-1}) 
  \stackrel{\mbox{\rm b}^+_2}{\rightarrow}
 (\sigma_1^{-1}\sigma_1^{-1}\sigma_1) 
 \stackrel{\mbox{\rm b}^-_2}{\rightarrow}
 (\sigma_1^{-1}\sigma_1) 
 \stackrel{\mbox{\rm II}^-_1}{\rightarrow} 
(\emptyset) 
]$
 $\stackrel{IIs}{\Longrightarrow}$
 
  \noindent
 $[(\emptyset) \stackrel{\mbox{\rm b}_0^+}{\rightarrow} 
 (\sigma_1^{-1}) \stackrel{\mbox{\rm b}_1^+}{\rightarrow}
 (\sigma_1^{-1}\sigma_1^{-1}) 
  \stackrel{\mbox{\rm b}^+_2}{\rightarrow}
 (\sigma_1^{-1}\sigma_1^{-1}\sigma_1)
  \stackrel{\mbox{\rm II}^-_2}{\rightarrow}
  (\sigma_1^{-1})
  \stackrel{\mbox{\rm II}^+_2}{\rightarrow}
 (\sigma_1^{-1}\sigma_1^{-1}\sigma_1)
  \stackrel{\mbox{\rm b}^-_2}{\rightarrow}
 (\sigma_1^{-1}\sigma_1) 
 \stackrel{\mbox{\rm II}^-_1}{\rightarrow} 
(\emptyset) 
]$
 $\stackrel{CC^2}{\Longrightarrow}$  
 
   \noindent
 $[(\emptyset) \stackrel{\mbox{\rm b}_0^+}{\rightarrow} 
 (\sigma_1^{-1}) \stackrel{\mbox{\rm b}_1^+}{\rightarrow}
 (\sigma_1^{-1}\sigma_1^{-1}) 
  \stackrel{\mbox{\rm b}^-_2}{\rightarrow}
 (\sigma_1^{-1})
  \stackrel{\mbox{\rm b}^+_1}{\rightarrow}
 (\sigma_1^{-1}\sigma_1) 
 \stackrel{\mbox{\rm II}_1^-}{\rightarrow} 
(\emptyset) 
]$
 $\stackrel{CC}{\Longrightarrow}$

   \noindent
 $[(\emptyset) \stackrel{\mbox{\rm b}_0^+}{\rightarrow} 
 (\sigma_1^{-1}) \stackrel{\mbox{\rm b}_1^+}{\rightarrow}
 (\sigma_1^{-1}\sigma_1^{-1}) 
  \stackrel{\mbox{\rm b}^-_2}{\rightarrow}
 (\sigma_1^{-1})
  \stackrel{\mbox{\rm b}^-_1}{\rightarrow}
(\emptyset) 
]$
 $\stackrel{2-H}{\Longrightarrow}$

    \noindent
 $[(\emptyset) \stackrel{\mbox{\rm b}_0^+}{\rightarrow} (\sigma_1^{-1})
  \stackrel{\mbox{\rm b}^-_1}{\rightarrow}
(\emptyset) 
]$
 $\stackrel{2-H}{\Longrightarrow}$
 
  \noindent
 $[
(\emptyset) 
].$
 
\begin{remark}{\rm
The regions between the attachment of the last $1$-handle and the first $2$-handle represent homeomorphisms of the torus. One can easily trace meridional and longitudinal classes through the homeomorphism and determine the lens space that is the $2$-fold branched cover of the knot $5_2$. In more generality, it is straightforward to determine the lens space structure of the $2$-fold branched cover of any $2$-bridge knot.

In even more generality, when we choose a braid representation of a knot, there is a standard Seifert surface associated to such a representative. A Heegaard splitting of the $2$-fold branched cover is determined since the optima of the diagram on the left of the diagram (recall we are arranging the knots horizontally) each corresponds to a $1$-handle attachment while the optima on the right correspond to $2$-handles. Attaching regions for $1$-handles encircle the arcs in the braid chart. The braid induces a homeomorphism of the Heegaard surface. And so a representation of the braid group in the mapping class group is determined.
}\end{remark}

\begin{remark}{\rm The Seifert surface need not be orientable in order to fold the $2$-fold branched cover. We can use unoriented charts, and map each covering surface in $\R^3$. The union will be a $3$-manifold mapped into $\R^4$. 
Furthermore, we can cut a non-orientable Seifert surface along a simple closed curve and obtain an orientable surface bounded by a link one component of which is the original knot. We use the orientation reversing loop to define a sequence of semi-oriented charts. 
Below the sequence of semi-oriented charts associated to the three-half twisted Mobius band whose boundary is the trefoil is indicated. 

Immediately following this remark, the braid essay that is associated is given. The crossing changes are indicated as $\sigma_1^{\pm} \stackrel{\chi^{-}_{j}}{\rightarrow} \sigma_1^{\mp}$.
The birth and death of such a curve are indicated by $\stackrel{\Xi_\pm}{\Longrightarrow}$ 
}\end{remark}

\noindent
{\bf The braid essay associated to the three-half twisted Mobius band:}

\vspace{.5cm}

\noindent 
\begin{wrapfigure}[11]{r}{3in}\vspace{-.4in}
\includegraphics[scale=.2]{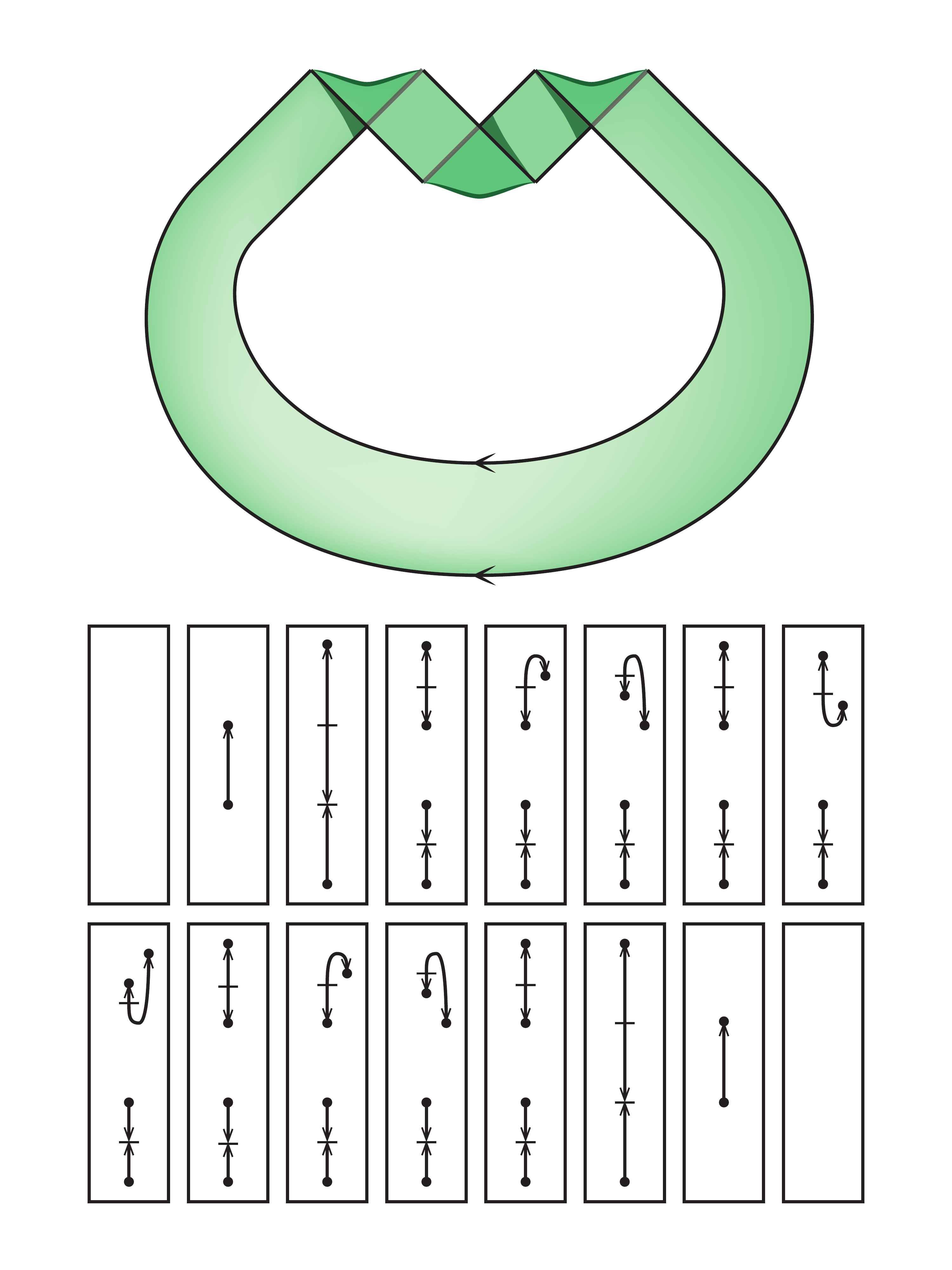}\vspace{-.4in}
\end{wrapfigure}
 $[(\emptyset)]$ $\stackrel{1-H}{\Longrightarrow}$
 
 \noindent
 $[(\emptyset) \stackrel{\mbox{\rm b}_0^+}{\rightarrow} 
 (\sigma_1) \stackrel{\mbox{\rm b}_1^-}{\rightarrow}
(\emptyset)]$
 $\stackrel{\Xi_+}{\Longrightarrow}$
 
 \noindent
 $[(\emptyset) \stackrel{\mbox{\rm b}_0^+}{\rightarrow} 
 (\sigma_1)  \stackrel{\chi_1^-}{\rightarrow} 
 (\sigma_1^{-1})   \stackrel{\chi_1^+}{\rightarrow}  (\sigma_1)
 \stackrel{\mbox{\rm b}_1^-}{\rightarrow}
(\emptyset)]$
 $\stackrel{1-H}{\Longrightarrow}$
 
 \noindent
 $[(\emptyset) \stackrel{\mbox{\rm b}_0^+}{\rightarrow} 
 (\sigma_1)  \stackrel{\chi_1^-}{\rightarrow} 
 (\sigma_1^{-1})    \stackrel{\mbox{\rm b}_1^-}{\rightarrow}
 (\emptyset) \stackrel{\mbox{\rm b}_0^+}{\rightarrow}
 (\sigma_1^{-1}) 
 \stackrel{\chi_1^+}{\rightarrow}  (\sigma_1)
 \stackrel{\mbox{\rm b}_1^-}{\rightarrow}
(\emptyset)]$
 $\stackrel{CC}{\Longrightarrow}$

 \noindent
 $[(\emptyset) \stackrel{\mbox{\rm II}_1^+}{\rightarrow} 
 (\sigma_1 \sigma_1^{-1})  \stackrel{\mbox{\rm b}_2^-}{\rightarrow} 
 (\sigma_1)
 \stackrel{\chi_1^-}{\rightarrow} 
 (\sigma_1^{-1})    \stackrel{\mbox{\rm b}_1^-}{\rightarrow} 
  (\emptyset) \stackrel{\mbox{\rm b}_0^+}{\rightarrow}
 (\sigma_1^{-1}) 
 \stackrel{\chi_1^+}{\rightarrow}  (\sigma_1)
 \stackrel{\mbox{\rm b}_1^-}{\rightarrow}
(\emptyset)]$
 $\stackrel{X}{\Longrightarrow}$

 \noindent
 $[(\emptyset) \stackrel{\mbox{\rm II}_1^+}{\rightarrow} 
(\sigma_1 \sigma_1^{-1}) 
 \stackrel{\chi_1^-}{\rightarrow} 
(\sigma_1^{-1} \sigma_1^{-1})
\stackrel{\mbox{\rm b}_1^-}{\rightarrow} 
 (\sigma_1^{-1})
   \stackrel{\mbox{\rm b}_1^-}{\rightarrow} 
  (\emptyset) \stackrel{\mbox{\rm b}_0^+}{\rightarrow}
 (\sigma_1^{-1}) 
 \stackrel{\chi_1^+}{\rightarrow}  (\sigma_1)
 \stackrel{\mbox{\rm b}_1^-}{\rightarrow}
(\emptyset)]$
 $\stackrel{CC}{\Longrightarrow}$
 
  \noindent
 $[(\emptyset) \stackrel{\mbox{\rm b}_0^+}{\rightarrow} 
 (\sigma_1)  \stackrel{\chi_1^-}{\rightarrow} 
 (\sigma_1^{-1})    \stackrel{\mbox{\rm b}_1^-}{\rightarrow}
 (\emptyset) \stackrel{\mbox{\rm b}_0^+}{\rightarrow}
 (\sigma_1^{-1}) 
 \stackrel{\chi_1^+}{\rightarrow}  (\sigma_1)
 \stackrel{\mbox{\rm b}_1^-}{\rightarrow}
(\emptyset)]$
 $\stackrel{CC}{\Longrightarrow}$

   \noindent
 $[(\emptyset) \stackrel{\mbox{\rm b}_0^+}{\rightarrow} 
 (\sigma_1)  \stackrel{\chi_1^-}{\rightarrow} 
 (\sigma_1^{-1})    \stackrel{\mbox{\rm b}_1^+}{\rightarrow} (\sigma_1^{-1} \sigma_1) \stackrel{\mbox{\rm II}_1^-}{\rightarrow} 
 (\emptyset) \stackrel{\mbox{\rm b}_0^+}{\rightarrow}
 (\sigma_1^{-1}) 
 \stackrel{\chi_1^+}{\rightarrow}  (\sigma_1)
 \stackrel{\mbox{\rm b}_1^-}{\rightarrow}
(\emptyset)]$
 $\stackrel{X}{\Longrightarrow}$

    \noindent
 $[(\emptyset) \stackrel{\mbox{\rm b}_0^+}{\rightarrow} 
 (\sigma_1) \stackrel{\mbox{\rm b}_0^+}{\rightarrow} 
 (\sigma_1 \sigma_1) 
  \stackrel{\chi_1^-}{\rightarrow} 
 (\sigma_1^{-1} \sigma_1)\stackrel{\mbox{\rm II}_1^-}{\rightarrow} 
 (\emptyset) \stackrel{\mbox{\rm b}_0^+}{\rightarrow}
 (\sigma_1^{-1}) 
 \stackrel{\chi_1^+}{\rightarrow}  (\sigma_1)
 \stackrel{\mbox{\rm b}_1^-}{\rightarrow}
(\emptyset)]$
 $\stackrel{CC}{\Longrightarrow}$

 \noindent
 $[(\emptyset) \stackrel{\mbox{\rm b}_0^+}{\rightarrow} 
 (\sigma_1)  \stackrel{\chi_1^-}{\rightarrow} 
 (\sigma_1^{-1})    \stackrel{\mbox{\rm b}_1^-}{\rightarrow}
 (\emptyset) \stackrel{\mbox{\rm b}_0^+}{\rightarrow}
 (\sigma_1^{-1}) 
 \stackrel{\chi_1^+}{\rightarrow}  (\sigma_1)
 \stackrel{\mbox{\rm b}_1^-}{\rightarrow}
(\emptyset)]$
 $\stackrel{CC}{\Longrightarrow}$

 \noindent
 $[(\emptyset) \stackrel{\mbox{\rm II}_1^+}{\rightarrow} 
 (\sigma_1 \sigma_1^{-1})  \stackrel{\mbox{\rm b}_2^-}{\rightarrow} 
 (\sigma_1)
 \stackrel{\chi_1^-}{\rightarrow} 
 (\sigma_1^{-1})    \stackrel{\mbox{\rm b}_1^-}{\rightarrow} 
  (\emptyset) \stackrel{\mbox{\rm b}_0^+}{\rightarrow}
 (\sigma_1^{-1}) 
 \stackrel{\chi_1^+}{\rightarrow}  (\sigma_1)
 \stackrel{\mbox{\rm b}_1^-}{\rightarrow}
(\emptyset)]$
 $\stackrel{X}{\Longrightarrow}$

 \noindent
 $[(\emptyset) \stackrel{\mbox{\rm II}_1^+}{\rightarrow} 
(\sigma_1 \sigma_1^{-1}) 
 \stackrel{\chi_1^-}{\rightarrow} 
(\sigma_1^{-1} \sigma_1^{-1})
\stackrel{\mbox{\rm b}_1^-}{\rightarrow} 
 (\sigma_1^{-1})
   \stackrel{\mbox{\rm b}_1^-}{\rightarrow} 
  (\emptyset) \stackrel{\mbox{\rm b}_0^+}{\rightarrow}
 (\sigma_1^{-1}) 
 \stackrel{\chi_1^+}{\rightarrow}  (\sigma_1)
 \stackrel{\mbox{\rm b}_1^-}{\rightarrow}
(\emptyset)]$
 $\stackrel{CC}{\Longrightarrow}$
 
  \noindent
 $[(\emptyset) \stackrel{\mbox{\rm b}_0^+}{\rightarrow} 
 (\sigma_1)  \stackrel{\chi_1^-}{\rightarrow} 
 (\sigma_1^{-1})    \stackrel{\mbox{\rm b}_1^-}{\rightarrow}
 (\emptyset) \stackrel{\mbox{\rm b}_0^+}{\rightarrow}
 (\sigma_1^{-1}) 
 \stackrel{\chi_1^+}{\rightarrow}  (\sigma_1)
 \stackrel{\mbox{\rm b}_1^-}{\rightarrow}
(\emptyset)]$
 $\stackrel{2-H}{\Longrightarrow}$

 \noindent
 $[(\emptyset) \stackrel{\mbox{\rm b}_0^+}{\rightarrow} 
 (\sigma_1)  \stackrel{\chi_1^-}{\rightarrow} 
 (\sigma_1^{-1})   \stackrel{\chi_1^+}{\rightarrow}  (\sigma_1)
 \stackrel{\mbox{\rm b}_1^-}{\rightarrow}
(\emptyset)]$
 $\stackrel{\Xi_-}{\Longrightarrow}$

 \noindent
 $[(\emptyset) \stackrel{\mbox{\rm b}_0^+}{\rightarrow} 
 (\sigma_1) \stackrel{\mbox{\rm b}_1^-}{\rightarrow}
(\emptyset)]$
 $\stackrel{2-H}{\Longrightarrow}$
 
 \noindent
 $[(\emptyset)]$ $\stackrel{1-H}{\Longrightarrow}$

\section{The $2$-fold branched cover of $S^4$ branched along a knotted or linked surface}

In~\cite{CS:Seif} a Seifert algorithm for knotted surfaces that project without branch points is presented. This algorithm was adjusted by the second author in~\cite{Kam:Seif} to be applied to the surface braid case. In our respective books, \cite{CS1998} or \cite{Kam2002} these algorithms are described in detail. Since the $2$-fold branched cover of $S^4$ branched along an orientable knotted or linked surface can be constructed from two copies of  $S^4$ both cut along the Seifert solid, and gluing the positive side of one copy to the negative side of the other. 

The braid chart can be used to give a movie description of the branch locus. In such a movie, or in a standard movie description of the knotted surface, the Seifert solid is cut into $2$-dimensional slices that are Seifert surfaces for the classical cross-sectional links. Thus for any such cross-section, the previous section provided an embedding of the $2$-fold branched cover of that cross-section. These are, then, connected by either handle attachments or isotopies that are induced by the Reidemeister moves. 

Specifically, the birth of a simple closed curve corresponds to a $1$-handle attached between the successive $3$-dimensional $2$-fold branched covers. A $1$-handle attached between successive movie stills corresponds to a $2$-handle attached between  successive $3$-dimensional $2$-fold branched covers.  

Following the discussion of the $3$-fold branched covers of classical knots, two examples of $3$-fold branched covers (one an immersed folding and one an embedded folding) will be given. These examples are more complicated than the $2$-fold branched coverings. 

Here we indicate that the $2$-fold branched cover of $S^4$ branched along an unknotted sphere is $S^4$, and we indicate an embedded folding. Meanwhile, the $2$-fold branched cover of $S^4$ branched along an unknotted torus is $S^2 \times S^2$ and also has an embedded folding.

\begin{center}
\includegraphics[width=3in]{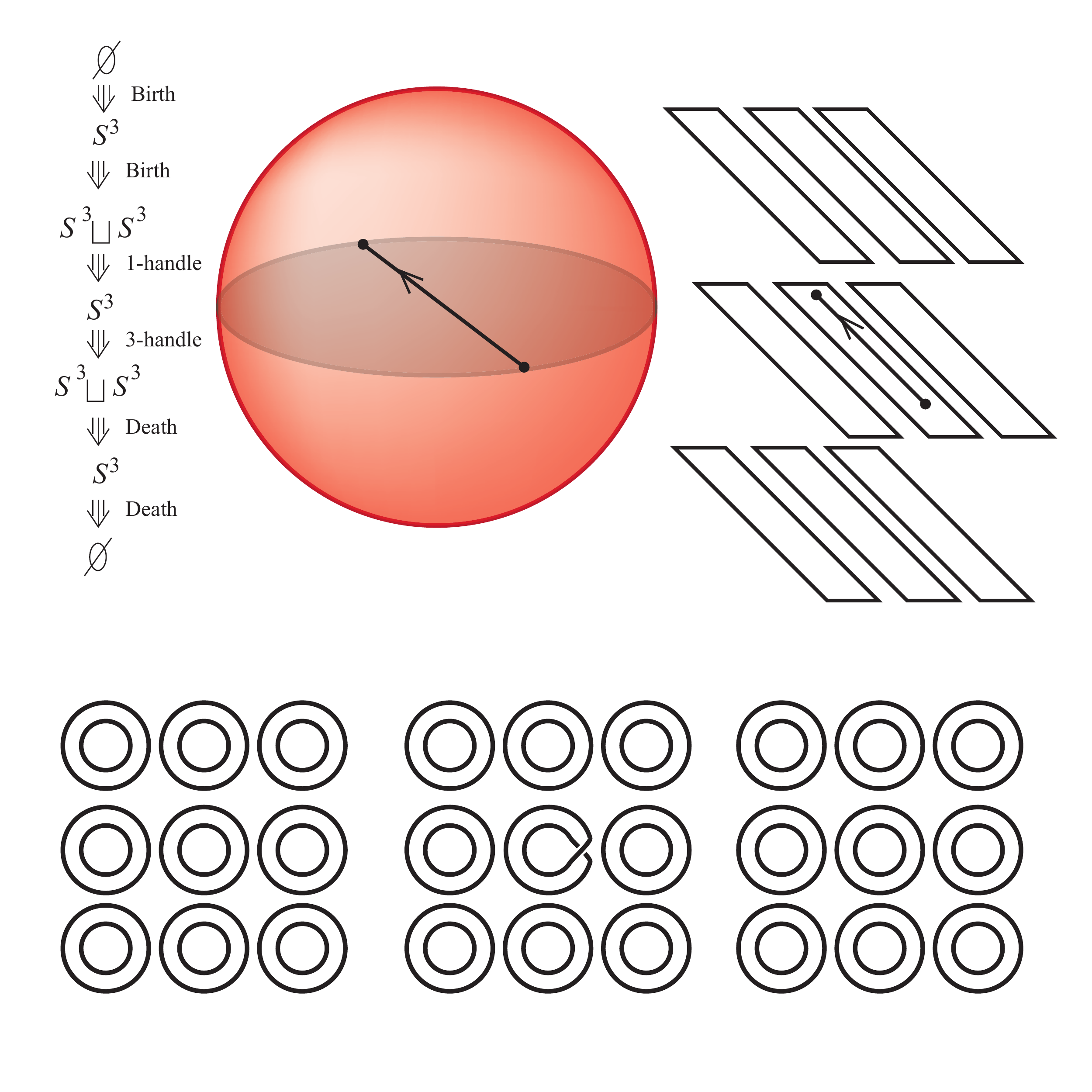}
\end{center}

 To summarize the constructions, the Seifert surface in dimension $3$ creates $2$-dimensional chart or curtain. By placing this in general position with respect to a height direction, we can cut between critical levels and connect the resulting sequence of $2$-dimensional charts via well-understood critical phenomena that correspond to chart moves. In this way, we can construct embedded foldings for the $2$-fold branched cover of $S^3$ branched along a knot or link. Furthermore, we can use non-orientable Seifert surfaces to construct immersed foldings. In $4$-space, we can use the Seifert algorithm for knotted surfaces to construct a Seifert solid. This solid can be expressed in movie form or via the surface braid picture of the knotted orientable surface. The critical events and Reidemeister moves have well-understood effects upon the cross-sectional Seifert surfaces. A folding of the $2$-fold branched cover of $S^4$ branched along a given surface knot or link is constructed, then, via interconnecting the corresponding covers of $S^3$. Thus we have proven Theorem~\ref{twofold}.

 \section{Three-fold simple branched covers}
 
 \begin{wrapfigure}[14]{l}{3.75in}\vspace{-.1in}
\includegraphics[scale=.09]{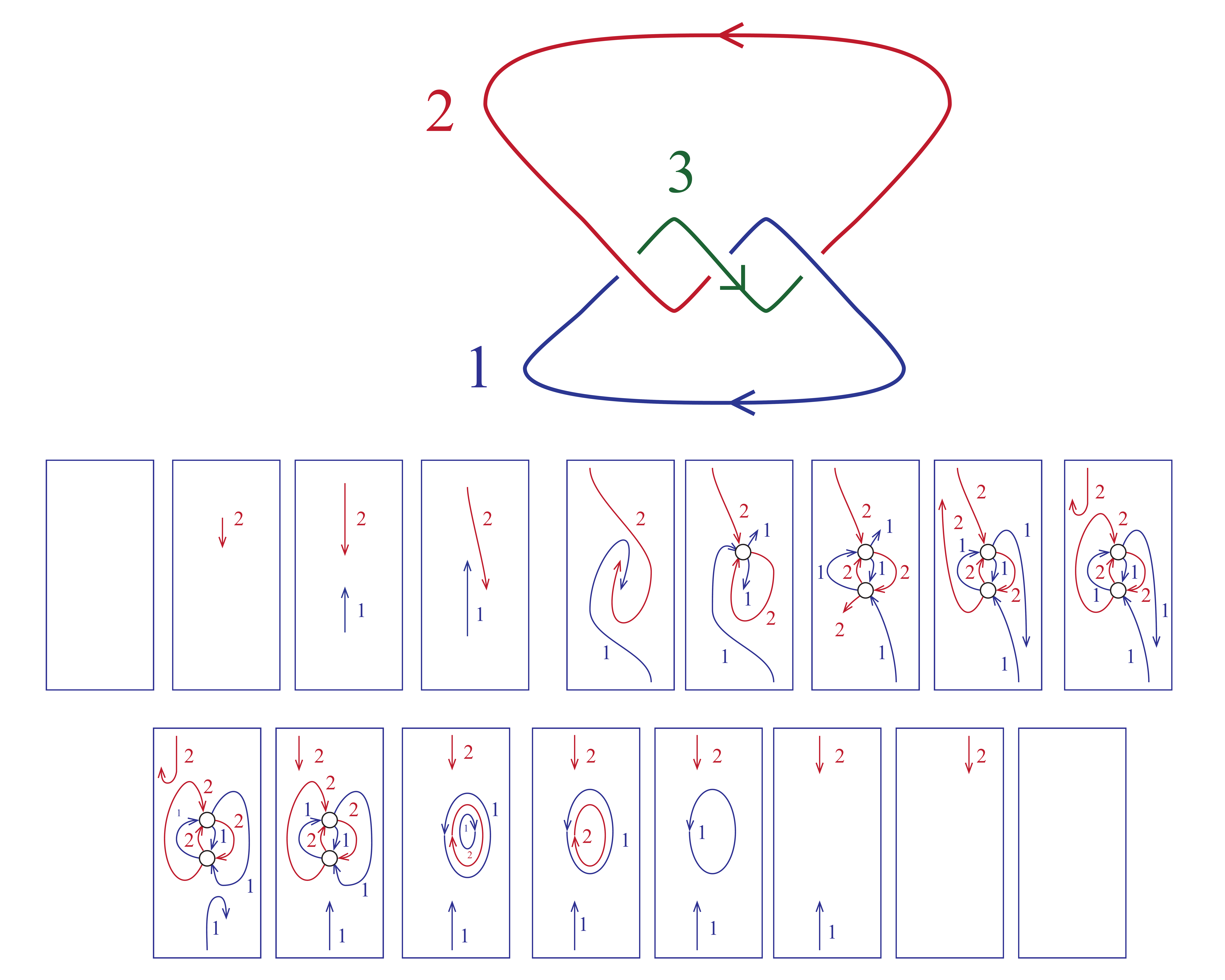}\end{wrapfigure}
 In this section, we construct embeddings and immersions of $3$-fold branch covers of the $3$ and $4$-sphere branched over the trefoil, the knot $7_4$, and the knotted sphere that is the spin of the trefoil. We begin with the trefoil as indicated to the left of this paragraph. The trefoil is $3$-colorable, and we can lift this coloring to a representation to the $3$-string braid group (also the fundamental group of the knot itself) by assigning the braid generator $\sigma_1$ to the (blue) label $1$, assigning the braid generator $\sigma_2$ to the (red) label $2$, and assigning $\sigma_2^{-1} \sigma_1 \sigma_2$ to the  (green) label $3$. From this coloring, a sequence of braid charts is constructed. The endpoints of the arcs in the charts correspond to black vertices, and as the charts are stacked, these end points trace the outline of the knot. Successive braid charts differ by a planar isotopy, by one of the chart moves, or by the addition or subtraction of a simple arc. As before, the addition or subtraction of an arc corresponds to attaching a handle between successive sheets in the covering. The sequence of charts can be interpolated in $3$-space to create a (non-generically) immersed surface in $3$-space that we call a {\it curtain}. To construct the $3$-fold branch cover, we take three copies of $S^3$, cut each along the appropriately labeled sheets  of the curtains, and re-glue. The embedding in this case is achieved by interpreting each chart as a braided surface and interpolating between successive surfaces.

  \begin{wrapfigure}{r}{3.5in}\vspace{-0.2in}
\includegraphics[scale=.052]{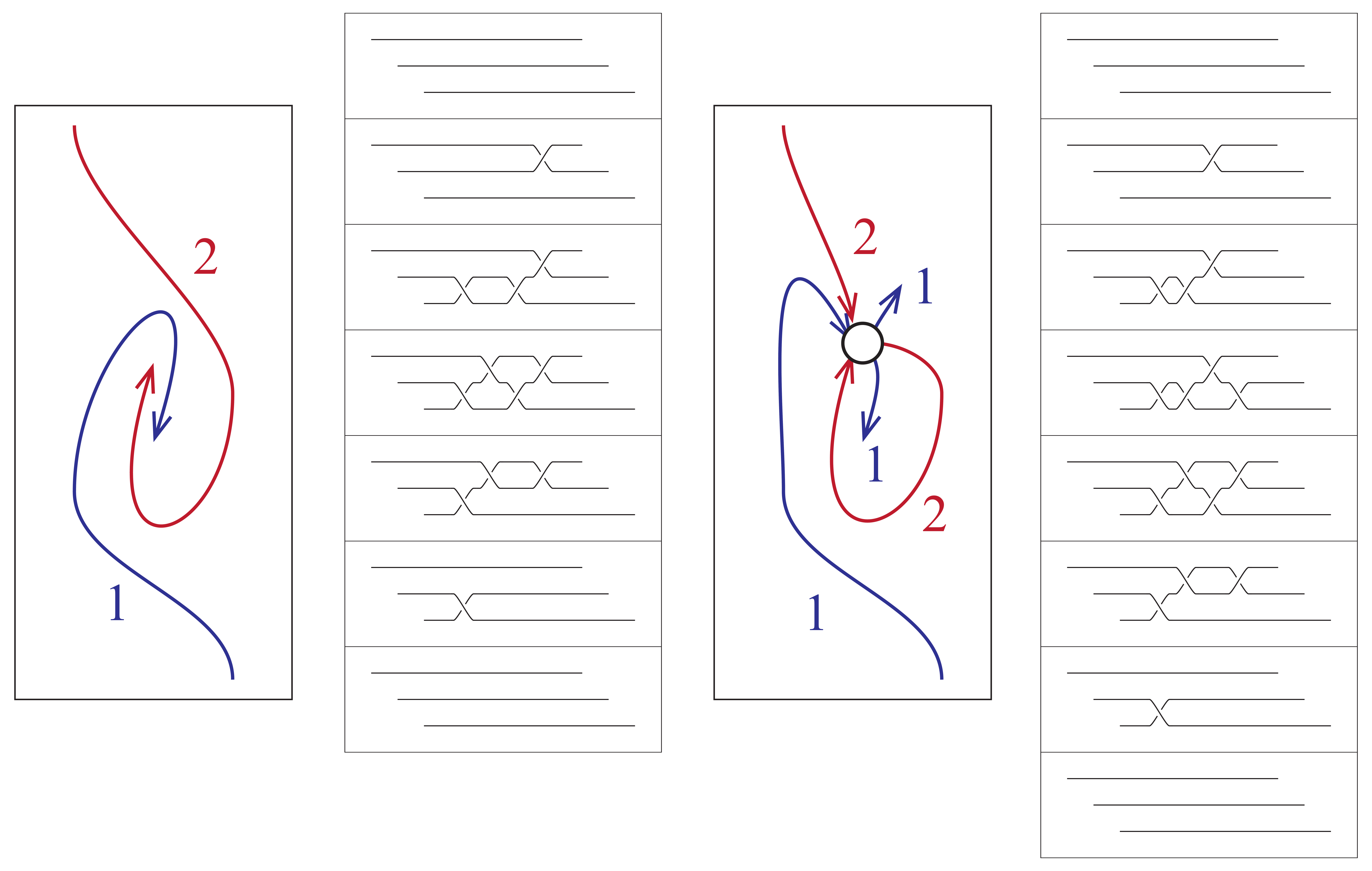}
\end{wrapfigure} For two cross-sections of the curtains, we have indicated the resulting braid movies. The difference between these two movies is a C-III move that corresponds to pushing a branch point through a transverse sheet. At the level of surface braids, this is a manifestation of one of the Roseman moves. The complete sequence of such surface isotopies and handle additions describes the branched cover explicitly embedded in $5$-space.

\begin{wrapfigure}[13]{l}{4in}\hspace{-0.15in}\vspace{-0.6in}\includegraphics[scale=0.09]{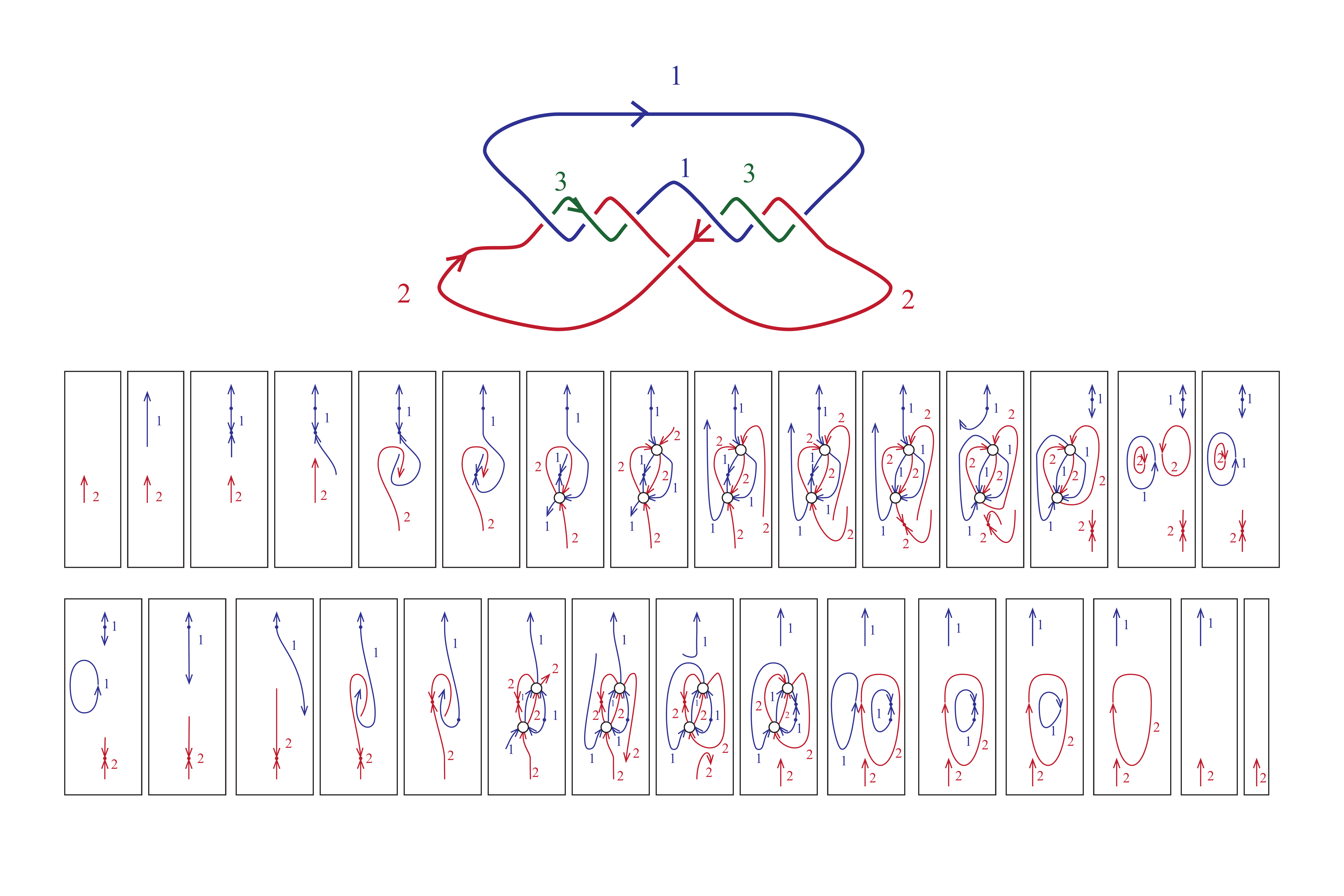}\end{wrapfigure}
To the left of this paragraph, the corresponding immersed folding of the $3$-fold branched cover of $S^3$ branched along the knot $7_4$ is illustrated. As before, the (blue) label $1$ indicates the braid generator $\sigma_1$ and the (red) label $2$ indicates $\sigma_2$. Here, however, nodes are introduced in the cross-sectional charts. Such a node indicates that a braid generator is switched with its inverse. In a cross-sectional surface braid, there is a self-transverse double point at such a node. In the case at hand there is a simple closed curve of such double points for the immersed folding of branched cover. It is known \cite{KitSuz2005}  that there is no surjective homomorphism between the fundamental groups of the complement of $7_4$ and the complement of the trefoil in which the peripheral structure is preserved. Such a homomorphism is necessary for there to be an embedded folding. See also \cite{HKMS2011}.

\begin{center}\includegraphics[scale=0.069]{spuntrefoil}\end{center}
Our last example for this paper is a construction of an embedded folding of the $3$-fold branch cover of $S^4$ branched along the spun trefoil. A movie for the spun-trefoil is palindromic with the square knot as the central step. We illustrate a curtain for each non-trivial cross-section. Each such curtain determines an embedding of the $3$-fold branched cover of $S^3$ branched along the corresponding knot or link (in fact trivial link until the middle level). Two curtains in successive cross-sections differ by easy to understand moves to curtains. For example, moving from the empty curtain to the curtain that is a red disk corresponds to attaching a $1$-handle between the second and third copy of $S^3$. The next few steps all correspond to performing a move and then undoing the same move. Thus the curtain move is to replace the identity with a pair of a move and its inverse. When a saddle is attached to the cross-sectional knot, a critical point of index $2$ is attached between the curtains.

\subsection{Proof of Theorem~\ref{mainimmersed} when $k=1$}

Suppose that a knot or link $K$ is given that is $3$-colorable. Assume that the knot is given in braid form with braid index $m$. By convention, we write the knot as the closure of the braid $\beta$ with the closing strings on the right of the braid, and orient the braid downward. 

The coloring induces a {\it color vector} at the top of the braid. This is the sequence $(c_1, \ldots, c_m)$ of colors at the top-left of the braid with $c_j \in \{1,2,3\}$. Each braid generator or its inverse induces a transformation of this vector. In a few sentences, we will describe these transformations in terms of the curtains. First, we use the color vector to describe the structure of the curtains before the braiding occurs. 

We may recolor the knot so that the first element, $c_1$, in the color vector is $1$. By convention we indicate this element by the color blue in our illustrations. If $c_j= 2$, then we color this red. The color $1$ corresponds to the transposition $(1,2)$ while $2$ corresponds to $(2,3)$. The remaining color $3$ corresponds to $(1,3)$; in the illustrations this is indicated in green. Observe that in the permutation group $(1,3)=(2,3)^{-1}(1,2)(2,3)$.

In anticipation of a lifting of the cross-sectional chart of a curtain to the braid group, we construct a specific chart that reflects the color vector immediately before braiding occurs. There is a sequence of horizontal arcs oriented as right-pointing arrows. Using a standard coordinate system, the top most arc points from $(-m,m)$ towards $(m,m)$. The arc immediately below this points from $(-m+1,m-1)$ towards $(m-1,m-1)$. Continue shortening the next arc by two units and move it one unit below the prior arc. Thus the lowest arc point from $(-1,1)$ to $(1,1)$. This top arc is labeled $1$ (colored blue).  Similarly, if $c_j=1$, then the arrow from $(-m-1+j,m+1- j)$ to $(m+1-j,m+1-j)$ is labeled $1$ (colored blue). If $c_j=2$, then the arrow from $(-m-1+j,m+1- j)$ to $(m+1-j,m+1-j)$ is labeled $2$ (colored red). Finally, if $c_j=3$, then the arrow from $(-m-1+j,m+1-j)$ to $(m+1-j,m+1-j)$ is colored $1$ but it is encircled by an anti-clockwise oriented oval that is colored $2$.

\begin{wrapfigure}[8]{l}{2.6in}\vspace{-0.2in}\includegraphics[scale=0.3]{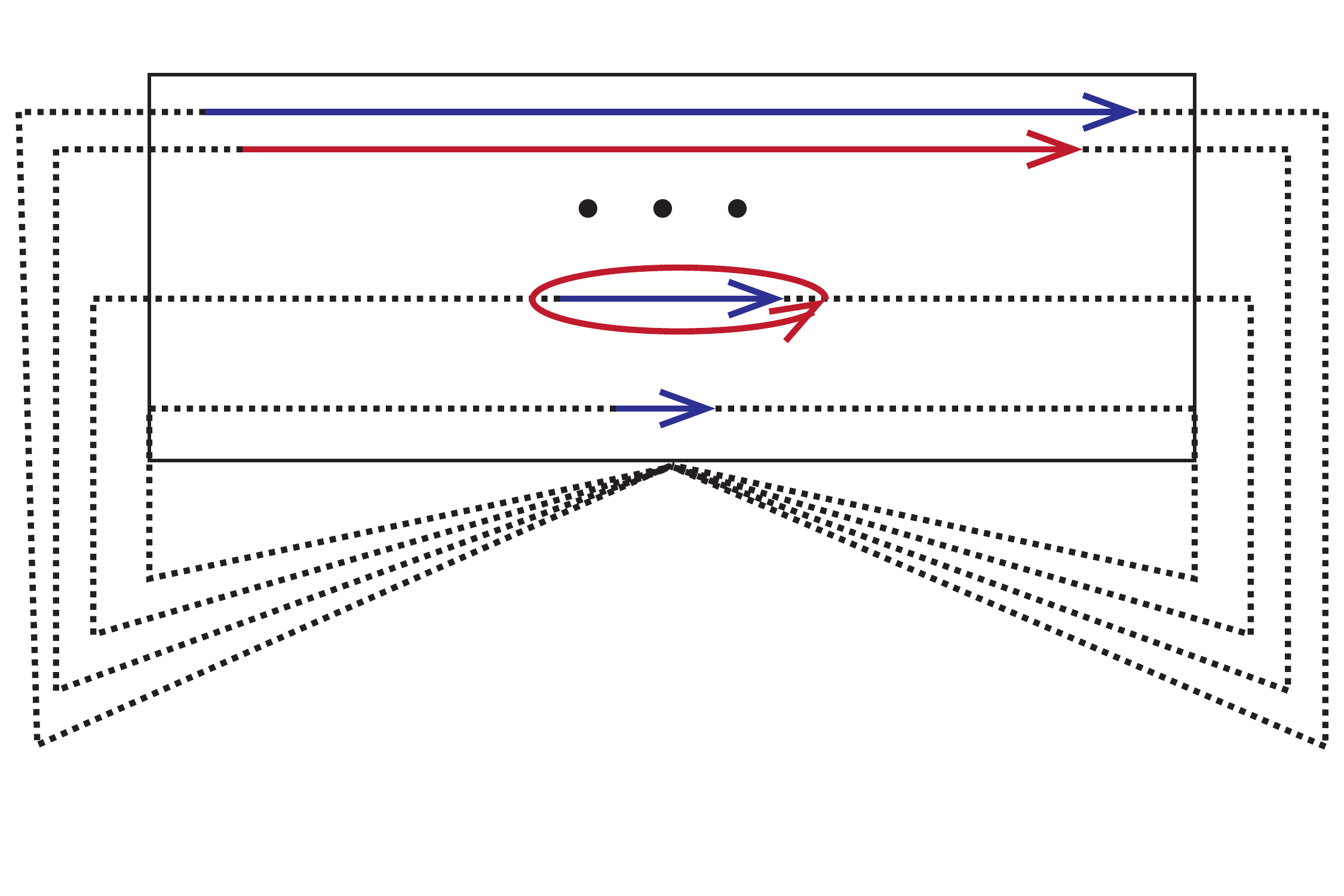}\end{wrapfigure} The figure to the left of this paragraph indicates the cross-section of the curtain before the braiding. Also indicated here by means of dotted arcs is an Hurwitz arc system that connects the end-points of the arrows to the base point which we take to be the origin. Moving up from this level, at the $j\/$th level the $j\/$th arc vanishes. If it is encircled by an oval, then this oval vanishes between the $j\/$th and $(j+1)\/$st level.  Thus before the knot appears, the curtain is empty and represents three nested $2$-spheres. After the first arc is born, the first and second spheres are connected by a $1$-handle. As each maximal point of the knot (in braid form) is passed, another $1$-handle is attached to the surface above. The ovals correspond to a type II, type II bubble move being performed so that a handle can be attached between the first and third sheets. The Hurwitz arcs on the left side of the figure record the color vector $(c_1, \ldots, c_m)$.

\begin{wrapfigure}{r}{4in}\includegraphics[scale=0.12]{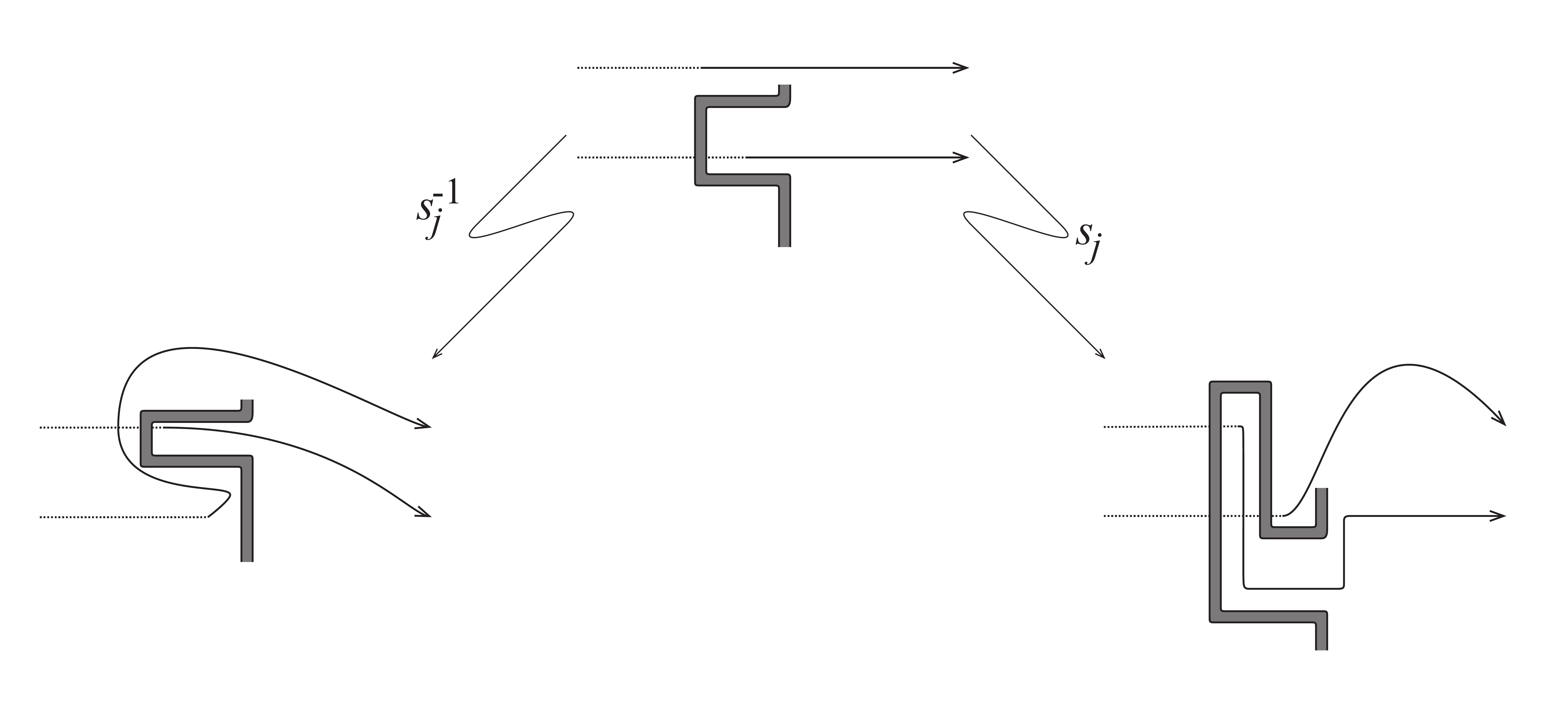}\end{wrapfigure} The figure to the right of this paragraph indicates the effect of the standard braid generator $s_j$ and its inverse to the curtains. The grey bands indicate that other arcs may cross the Hurwitz arc that connects to the black vertex. The result of the action of the braid generator or its inverse upon the color vector is determined by reading the intersection sequence with the corresponding Hurwitz arcs. 

At the bottom of the braid $\beta$, the curtain still consists of $k$ embedded arrows with some of them encircled, but the topography is potentially complicated. However, the color vector has returned to its original state $(c_1, \ldots, c_m)$. We proceed to simplify the curtain by applying chart moves, or, if necessary, by introducing nodes. Here is how to proceed.

\begin{wrapfigure}{l}{4in}\includegraphics[scale=0.27]{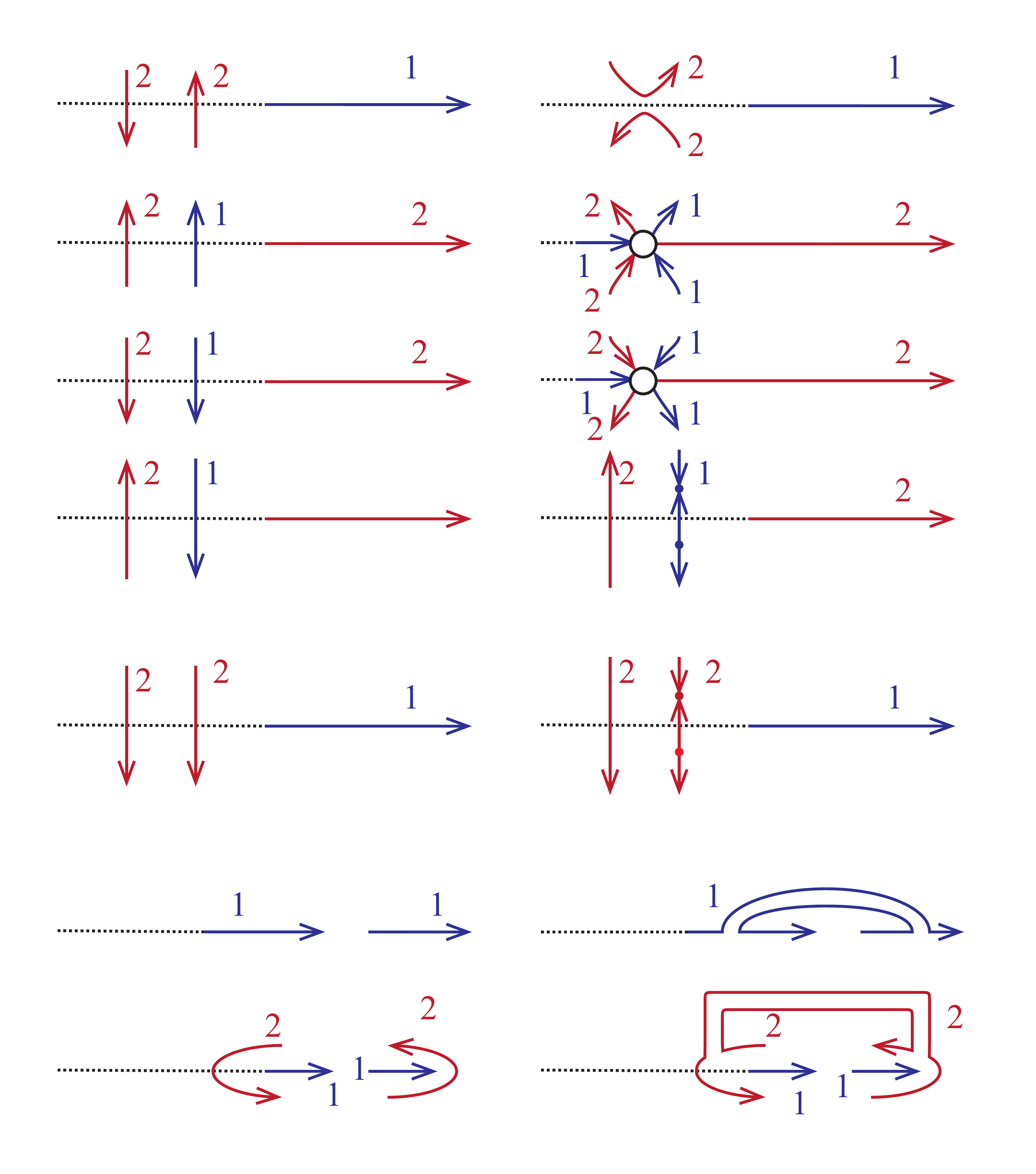}\end{wrapfigure} Consider the Hurwitz arc $\alpha_1$ that appears on the top left. It intersects the arcs and ovals of the chart in a sequence $w_1$. In this way the word $w_1 c_1 w_1^{-1}$ represents the element $c_1=(1,2)$ in the permutation group. Consequently, the intersection sequence $w_1$ represents a trivial word or $(1,2)$ in the permutation group $\Sigma_3$. The word $w_m$ also represents a word in the braid group since the orientations of the oval and the original arrows can be used to determine a direction of intersection with the Hurwitz arc. That is to say, at this point each chart above the bottom of the braid $\beta$ is oriented. We continue to consider the oriented case. If $w_1$ contains a syllable of the form $v v^{-1}$, where $v$ is red or blue, then we may eliminate this syllable by performing a saddle move and using a segment of the Hurwitz arc $\alpha_1$ as the core of the saddle. Such moves can be performed successively. Similarly, if $w_1$ ends in (red,blue) or in (red${}^{-1}$,blue${}^{-1}$), then a C-III move can be used to pull the black vertex at the end of the top arc to the left; doing so reduces the length of $w_1$. 

Continuing in this way, we gradually shorten $w_1$. However, $w_1$ may represent a trivial word in the permutation group or $(1,2)$, but may not represent a trivial word or $\sigma_1$ in the braid group. In that case, there is a syllable of the form $vv$ in the word $w$, and we introduce a pair of nodes on one of the segments and on either side of the Hurwitz arc $\alpha_1$. Doing so, changes the local orientation and the intersection sequence becomes $v^{-1}v$, so that a saddle cancelation can occur.  Thus we simplify the intersection word until it is empty or $\sigma_1$.  If the intersection word is $\sigma_1$, then by a saddle move, we may remove the intersection. At this stage, another saddle band can be attached, so that the points $(-m,m)$ and $(m,m)$ are joined by an embedded arc that does not intersect the rest of the chart. 

We proceed to successively simplify the intersection words between the Hurwitz arcs $\alpha_2$ through $\alpha_m$ so that they are each trivial. When each intersection sequence is empty, the left end of any arc or oval can be attached to the right end by means of saddle bands. We obtain a copy of the original chart immediately before the braiding and a disjoint closed chart that has nodes upon it. The closed chart can be eliminated by C-I moves, by moving nodes through white vertices (if necessary) and by canceling pairs of nodes. The copy of the original chart is trivialized by mimicking the moves that occur before the braid $\beta$. This completes the proof.

\section*{Future Work} In a subsequent work, we will demonstrate how to construct immersed foldings of the $3$-fold branched covers of $S^3$ branched along a knot or link. We are also examining simple branched covers of higher degree. In particular, we are interested in $5$-fold branch covers of $S^4$ in the light of Iori and Piergallini's theorem. We have constructed an immersion of the $3$-fold branched cover of $S^4$ branched along the $2$-twist-spun trefoil. These ideas lead naturally to a study of higher dimensional knotting via a braid theory. In particular, diagrammatic methods can be extended by using charts, curtains, and additional colorings.

\medskip

\begin{flushleft}
J. Scott Carter \\ 
Department of Mathematics \\ 
University of South Alabama \\ 
Mobile, AL 36688 \\
USA\\ 
E-mail address: {\tt carter@southalabama.edu} 
\end{flushleft}

\begin{flushleft}
Seiichi Kamada \\ 
Department of Mathematics \\ 
Hiroshima University  \\
Hiroshima 739-8526\\ 
JAPAN \\
E-mail address: {\tt  kamada@math.sci.hiroshima-u.ac.jp}
\end{flushleft}

\end{document}